\documentclass[lettersize,journal]{IEEEtran}

\usepackage{color,array}

\usepackage{graphicx}
\usepackage{amsmath,amsfonts}
\usepackage{amsthm}
\usepackage{algorithmic}
\usepackage{array}
\usepackage{textcomp}
\usepackage{stfloats}
\usepackage{url}
\usepackage{cite}
\usepackage{verbatim}
\usepackage{xcolor}

\usepackage{subcaption}

\usepackage{algorithm}

\allowdisplaybreaks[3]

\newcommand{\calG}{{\mathcal{G}}}

\newcommand{\R}{\mathbb{R}}
\newcommand{\N}{\mathcal{N}}

\newtheorem{theorem}{Theorem}[section]
\newtheorem{lemma}[theorem]{Lemma}

\newtheorem{proposition}[theorem]{Proposition}
\newtheorem{problem}{Problem}
\newtheorem{remark}[theorem]{Remark}
\newtheorem{assumption}{Assumption}

\newcommand{\complex}{\mathbb{C}}
\newcommand{\tsp}{\mathsf{T}}

\newcommand{\veci}{\boldsymbol{i}}
\newcommand{\vecv}{\boldsymbol{v}}
\newcommand{\vecV}{\boldsymbol{\nu}}
\newcommand{\vecI}{\boldsymbol{\iota}}

\newcommand{\vectheta}{\boldsymbol{\theta}}
\newcommand{\vecxi}{\boldsymbol{\xi}}
\newcommand{\vecs}{\boldsymbol{s}}
\newcommand{\vecy}{\boldsymbol{y}}

\newcommand{\vecp}{\boldsymbol{p}}
\newcommand{\vecq}{\boldsymbol{q}}
\newcommand{\vecpl}{\boldsymbol{p}_l}
\newcommand{\vecql}{\boldsymbol{q}_l}
\newcommand{\smally}{\boldsymbol{\bar{y}}}
\newcommand{\capitalY}{\boldsymbol{Y}}

\newcommand{\vecu}{\boldsymbol{u}}

\newcommand{\vecn}{\boldsymbol{n}}

\newcommand{\vecz}{\boldsymbol{z}}

\DeclareMathAlphabet{\mymathbb}{U}{BOONDOX-ds}{m}{n}

\begin{document}


\title{Learning to Pursue AC Optimal Power Flow Solutions with Feasibility Guarantees}
\author{Damola Ajeyemi, Yiting Chen, Antonin Colot, Jorge Cort\'{e}s, and Emiliano Dall'Anese 
\thanks{Damola Ajeyemi is with the Division of Systems Engineering, Boston University, Boston, MA 02215, USA (email: dajeyemi@bu.edu).}
\thanks{Yiting Chen is with the Department of Electrical and Computer Engineering, Boston University, Boston, MA 02215, USA (email: yich4684@bu.edu).}
\thanks{Antonin Colot is with the  University of Li\`{e}ge, B-4000 Liège, Belgium (email: antonin.colot@eliagrid-int.com).}
\thanks{Jorge Cort\'{e}s is with the Department of Mechanical and Aerospace Engineering, University of California San Diego, CA 92093 San Diego, USA (email: cortes@ucsd.edu).}
\thanks{Emilliano Dall'Anese is with the Department of Electrical and Computer Engineering and the Division of Systems Engineering, Boston University, Boston, MA 02215, USA (email: edallane@bu.edu).}
}

\maketitle

\begin{abstract}
This paper focuses on an AC optimal power flow (OPF) problem for distribution feeders equipped with controllable distributed energy resources (DERs). We consider a solution method that is based on a continuous approximation of the projected gradient flow -- referred to as the safe gradient flow -- that incorporates voltage and current information obtained either through real-time measurements or power flow computations. These two setups enable both feedback-based online and offline implementations. The safe gradient flow involves the solution of convex quadratic programs (QPs). To enhance computational efficiency, we propose a novel framework that employs a neural network approximation of the optimal solution map of the QP. The resulting  method has two key features: (a) it ensures that the DERs' setpoints are practically feasible, even for an online implementation or when an offline algorithm has an early termination; (b) it ensures convergence to a neighborhood of a strict local optimizer of the AC OPF. The proposed method is tested 
on a 93-node distribution system with realistic loads and renewable generation. The test shows that our method successfully regulates voltages within limits during periods with high renewable generation.
\end{abstract}


\section{Introduction}
\label{sec:introduction}

This work considers power distribution systems with controllable distributed energy resources (DERs), and aims to advance real-time control strategies and computational methodologies in this domain.  The focus is on the AC optimal power flow (OPF) problem~\cite{molzahn2017survey} and, in particular, on its real-time implementation. These include recent frameworks that leverage feedback-based  implementations~\cite{gan2016online,dall2016optimal,bernstein2019real,picallo2022adaptive} or low-latency batch solutions~\cite{venzke2020learning,zamzam2020learning}. These real-time implementations seek to generate setpoints at a time scale that is consistent with the variability of uncontrollable loads and power available from renewable sources~\cite{taylor2016power}.

\emph{Prior work}. Feedback-based online algorithms have been explored in the context of AC OPF for distribution systems~\cite{gan2016online,dall2016optimal,bernstein2019real,picallo2022adaptive}. Shifting from a feedback optimization paradigm to feedforward optimization, a substantial body of work has explored the use of neural networks and deep learning techniques to approximate solutions to the AC OPF problem; see, for example,~\cite{zamzam2020learning,singh2021learning,nellikkath2022physics,pan2022deepopf,pan2023deepopf,park2023compact,fioretto2020predicting,tran2024learning,baker2020learning,chen2025physics}, the generative model in~\cite{wang2022fast}, and the foundation models in~\cite{hamann2024foundation}. While these methods primarily target AC OPF tasks in transmission networks, some of them can be adapted to distribution grids as well.  This body of literature has adopted various approaches: some aim to directly predict a solution to the AC OPF problem~\cite{zamzam2020learning}, while others focus on predicting a Karush–Kuhn–Tucker (KKT) point~\cite{fioretto2020predicting}. 

In general, these methods lack formal guarantees in terms of generating optimal solutions of the AC OPF, not to mention feasible points (as we will show in our numerical results). Once a candidate solution is generated by the neural network, recovering a valid operating point that satisfies all AC OPF constraints can be computationally demanding 
--  offsetting the speed advantages offered by the neural network approximation; heuristics may be used, but still lack formal guarantees. Post-processing for neural networks approximating solutions to problems with linear constraints are available in the literature~\cite{li2023learning}, but they are not applicable to the AC OPF. 
The AC OPF is nonconvex and may admit multiple globally  and locally optimal solutions; this means that the function that maps loads in the network and parameters of the problem into optimal solutions is a set-valued mapping. Such a set-valued mapping cannot be approximated with the single-valued mapping of a neural network; see the discussion in, e.g.~\cite{sun2018learning} and~\cite{pan2023deepopf,zhou2020optimal}. One workaround suggested in~\cite{sun2018learning}  is when the number of solutions (or KKT points) of the AC OPF is finite and they can all be identified; in this case,  the neural network can be trained to output the vector enumerating the optimal solutions (or KKT points). Enumerating the solutions (or KKT points) of the AC OPF is computationally infeasible~\cite{zhou2020optimal}. 

An alternative strategy involves replacing algorithmic updates in traditional optimization methods -- such as Newton-type or gradient-based methods -- with neural networks~\cite{baker2020learning,chen2025physics}. These methods offer computational advantages, but existing works do not offer convergence and feasibility guarantees. 

\emph{Contributions}. In this paper, we consider a solution method for the AC OPF that is based on
the so called safe gradient flow~\cite{allibhoy2023control,colot2024optimal}, {a constrained gradient flow whose design employs Control Barrier Functions (CBFs),} modified to incorporate
voltage and current information obtained either through real-time measurements (as in feedback optimization~\cite{gan2016online,dall2016optimal,bernstein2019real,picallo2022adaptive}) or power flow computations. To favor computational efficiency and speed, we propose a novel framework that employs a neural network approximation of the safe gradient flow. In particular, the neural network predicts the unique optimal solution of a quadratic program (QP) defining the map of the safe gradient flow. The learning task is well-posed, in the sense that the optimal solution map of the QP is a single-valued function and it is continuous. The learned  safe gradient flow is then used in conjunction with voltage and current information to identify AC OPF solutions. We summarize our contributions as follows:

\emph{(c1)} We propose an iterative method where the neural network approximation of the  safe gradient
flow is used with either real-time measurements or power flow computations.  When using real-time measurements, the proposed method constitutes a neural network-based feedback optimization method.  

\emph{(c2)} We show that our method leads to solutions that are practically feasible.  The term practical  feasibility refers to the fact that we provide guarantees on  the maximum constraint violation; the analytical estimate of the violation  allows for a careful tightening of the constraints in the AC OPF so that the neural network can be trained not to violate the actual constraints. The practical feasibility is at \emph{any time}, in the sense that the algorithm produces feasible points even when terminated before convergence or implemented online. To the best of our knowledge, this is the first neural network-based algorithm for the AC OPF problem with feasibility guarantees, without postprocessing. 

\emph{(c3)} We show that the proposed learning-based method converges exponentially fast within a neighborhood of KKT points of the AC OPF that are strict local optimizers.

\emph{(c4)} Numerical simulations are performed on a 93-bus distribution grid \cite{SimBench}  with realistic time-series for load and renewable generation
from the Open Power System 
Data. The considered network \cite{SimBench} is a low voltage benchmark distribution system, reflecting realistic operational settings. We show that our approach ensures voltage regulation and satisfaction of the DERs' constraints. Our proposed method shows superior performance in terms of voltage
regulation compared to approaches that attempt to approximate the solutions of the OPF directly, even  when using more data to train the latter. 

Finally, we note that this paper provides several innovations relative to our previous work~\cite{colot2024optimal}. We propose a novel framework using neural network approximations of the safe gradient flow to solve the AC optimal power flow problem; the proposed method leads to implementations for both online optimization and offline computation, providing new computationally-efficient methods for both applications (as shown shortly in Figure~\ref{fig:proposed_arch}). We provide  results for convergence and feasibility that  are new in the literature, and we test the neural network-based solution on a voltage regulation problem. We provide tighter bounds when we characterize the forward invariance of an inflated set of feasible points. These improved results can be used to effectively train the proposed NN-SGF under tightened constraints to ensure forward invariance of the original feasible set. 

The remainder of the paper is organized as follows. Section~\ref{sec:problem-formulation} formulates the AC OPF and explains our proposed mathematical model. Section~\ref{sec:NNopf} provides details on the neural network-based safe gradient flow, while Section~\ref{sec:results} illustrates simulation results. Section~\ref{sec:convergence} presents our theoretical results and Section~\ref{sec:conclusions} gathers our conclusions. Finally, the appendix contains the technical proofs.\footnote{\emph{Notation}. We use the following notational conventions throughout the paper. Boldface upper-case letters (e.g., \(\mathbf{X}\)) denote matrices, and boldface lower-case letters (e.g., \(\mathbf{x}\)) denote column vectors. The transpose of a vector or matrix is denoted by \((\cdot)^\top\), and the complex conjugate by \((\cdot)^*\). The imaginary unit is denoted by \(j\), satisfying \(j^2 = -1\), and the absolute value of a scalar is written as \(|\cdot|\). For notational simplicity, given a vector-valued function $\mathbf{f}(\mathbf{x}): \mathbb{R}^n \rightarrow \mathbb{R}^m$, the notation $|\mathbf{f}(\mathbf{x})|$ means that the absolute value is taken entry-wise; i.e., $|\mathbf{f}(\mathbf{x})| = [|f_1(\mathbf{x})|, \ldots, |f_m(\mathbf{x})|]^\top$. For a real-valued vector \(\mathbf{x} \in \mathbb{R}^N\), \(\mathrm{diag}(\mathbf{x})\) returns an \(N \times N\) diagonal matrix with the entries of \(\mathbf{x}\) on the diagonal. The \(\ell_2\)-norm of a vector \(\mathbf{x} \in \mathbb{R}^n\) is denoted \(\|\mathbf{x}\|\); for a matrix $\mathbf{X} \in \mathbb{R}^{n \times m}$, $\|\mathbf{X}\|$ is the induced  \(\ell_2\)-norm. The \(\ell_1\)-norm of a vector \(\mathbf{x} \in \mathbb{R}^n\) is denoted \(\|\mathbf{x}\|_1\), and the \(\ell_\infty\)-norm  is denoted \(\|\mathbf{x}\|_\infty\). For two vectors \(\mathbf{x} \in \mathbb{R}^n\) and \(\mathbf{u} \in \mathbb{R}^m\), the notation \((\mathbf{x}, \mathbf{u}) \in \mathbb{R}^{n+m}\) denotes their concatenation. The symbol \(\mathbf{0}\) is used to denote vectors or matrices of zeros, with dimension determined from context. The set of complex numbers is denoted \(\mathbb{C}\). For a complex vector $\boldsymbol{x} \in \mathbb{C}^{N}$, $\Re(\boldsymbol{x}) \in \mathbb{R}^{N}$ denotes its real part and $\Im(\boldsymbol{x}) \in \mathbb{R}^{N}$ its imaginary part. We denote by $\mathbb{N}_0$ the set of non-negative integers, and by $\mathbb{N}_{>0}$ the set of positive integers. The set of all integers is denoted by \(\mathbb{Z}\). }

\section{Problem Formulation and Proposed Model}
\label{sec:problem-formulation}


The paper considers a distribution system with \(N+1\) nodes, labeled by \(\{0,1,\dots,N\}\). Node \(0\) represents the substation (or point of common coupling), whereas \(\N := \{1,\dots,N\}\) contains the remaining nodes; these nodes may feature a mix of uncontrollable loads and controllable DERs. We focus on a steady-state representation in which currents and voltages are modeled as complex phasors. For each node \(k \in \N\), let the line-to-ground voltage phasor be $v_k \in \mathbb{C}$. The current injected at node \(k\) is $i_k$. At the substation, the voltage is denoted as \(v_0 = V_0 e^{j\,\delta_0}\)~\cite{Kerstingbook}.
As usual, applying Ohm’s and Kirchhoff’s Laws in the phasor domain yields the relationship
\begin{equation}
\label{eqn:Ymatrix}
\begin{bmatrix}
i_0 \\[3pt]
\veci_\mathcal{N}
\end{bmatrix}
\;=\;
\begin{bmatrix}
y_0 & \smally^\top \\[3pt]
\smally & \capitalY
\end{bmatrix}
\begin{bmatrix}
v_0\\[3pt]
\vecv_\mathcal{N}
\end{bmatrix},
\end{equation}
where \(\veci_\mathcal{N}=[\,i_1,\dots,i_N\,]^\tsp\!\in\complex^N\) and \(\vecv_\mathcal{N}=[\,v_1,\dots,v_N\,]^\tsp\!\in\complex^N\), and where the admittance matrix \(\capitalY\!\in\!\complex^{N\times N}\) and the vectors \(\smally\!\in\!\complex^N,y_0\in\complex\) are built based on the series and shunt parameters of the lines under a \(\Pi\)-model~\cite{Kerstingbook}. 
 
Suppose that there are \(G\) DERs in the network, each capable of generating or consuming active and reactive powers. Let $\vecu \;=\; [\,p_1,\dots,p_G,\;q_1,\dots,q_G\,]^\tsp\;\in\;\R^{2G}$ collect the DERs' active powers \(p_i\) and reactive powers \(q_i\). For each DER \( i \in \mathcal{G} \), the set of admissible active and reactive power setpoints is defined by a compact set \( \mathcal{C}_i \subset \mathbb{R}^2 \); the overall control domain is given by the Cartesian product \( \mathcal{C} := \mathcal{C}_1 \times \mathcal{C}_2 \times \dots \times \mathcal{C}_G \subset \mathbb{R}^{2G} \). Define the mapping \(m:\{1,\dots,G\}\to\N\) to indicate the node at which each DER is connected. Then, the net injections at node \(n\) can be written as $p_{\mathrm{net},n}
=\;
\sum_{i\in\calG_{n}}p_i
\;-\;p_{\ell,n}$, $q_{\mathrm{net},n}
=\;
\sum_{i\in\calG_{n}}q_i
\;-\;q_{\ell,n}$,
where \(\calG_n=\{\,i\in\{1,\dots,G\}:m(i)=n\}\) and with $p_{\ell,n}, q_{\ell,n}$ denoting the real and reactive loads (positive entries imply consumption).
Let $\vecp \in \mathbb{R}^N$ and $\vecq \in \mathbb{R}^N$ collect the active and reactive powers from the DERs on nodes $n \in \mathcal{N}$. Then, from~\eqref{eqn:Ymatrix}, one can derive the equation:
\begin{equation}
\label{eqn:power-flow}
(\vecp - \vecp_l) + j (\vecq - \vecq_l) \;=\;\mathrm{diag}(\vecv)\,\bigl(\smally^{*} v_0^{*} \;+\;\capitalY^*\,\vecv^{*}\bigr),
\end{equation}
where $\vecp_l, \vecq_l \in \mathbb{R}^N$ are the vectors collecting the aggregate active and reactive powers, respectively, of non-controllable loads at each node. Hereafter, we let $\vecs_l := [\vecp_l^\top, \vecq_l^\top]^\top$. Finally, nodes $\mathcal{M} \subseteq \mathcal{N}$, $M = |\mathcal{M}|$, are monitored and voltages can be measured, and lines $\mathcal{L} = \{1, \ldots, L\} \subset \mathcal{N} \times \mathcal{N}$ are also monitored and currents can be measured.

{In this formulation, we consider  a standard setup where node $0$ is designated as the  reference bus with $V_0=1$ p.u. and $\delta_0=0$. All remaining nodes $k\in\mathcal{N}$ are treated as PQ buses. Inverter-interfaced distributed energy resources (DERs) are represented as controllable active and reactive power injections $(p_i,q_i)$  at their respective interconnection nodes. Power injections are defined on a net basis, accounting for all elements connected at each node: injections are positive at nodes hosting DERs, negative at pure load nodes, and net values at nodes hosting both generation and demand. 
} Given the controllable powers $\vecu$ and the loads $\vecp_l, \vecq_l$, one can employ numerical techniques
to solve~\eqref{eqn:power-flow} for the voltages $\vecv$ {(for example, the Newton--Raphson algorithm or fixed-point methods)~\cite{chen2023fixed,wang2017existence}}. It is important to note that the power flow equation~\eqref{eqn:power-flow} may admit zero, one, or multiple solutions~\cite{bolognani2015existence,bernstein2018load,wang2017existence}. If multiple solutions exist, we focus on \emph{practical} solutions; i.e., the solution within the neighborhood of the nominal voltage profile that yields relatively high voltage magnitudes and low line currents. Due to the Implicit Function Theorem, we can define a map $(\vecu, \vecs_l) \mapsto \vecv(\vecu, \vecs_l) := [v_1(\vecu, \vecs_l), \ldots, v_{M}(\vecu, \vecs_l)]^\top$, mapping loads and power from the DERs into voltages at the monitored nodes. Additionally, based on the function $\vecv(\vecu, \vecs_l)$ and the topology of the network, we also define the function $(\vecu, \vecs_l) \mapsto \veci(\vecu, \vecs_l) := [i_1(\vecu, \vecs_l), \ldots, i_{L}(\vecu, \vecs_l)]^\top$, mapping loads and DERs' powers into monitored line currents.

The functions $\vecv(\vecu, \vecs_l)$ and $\veci(\vecu, \vecs_l)$ are utilized to formulate instances of the AC OPF. We adopt the notation $|\vecv(\vecu,\vecs_l)|$ and $|\veci(\vecu,\vecs_l)|$ to denote the entry-wise absolute value of these vector-valued functions; i.e., $|\vecv(\vecu,\vecs_l)| = [|v_1(\vecu,\vecs_l)|, \dots, |v_M(\vecu,\vecs_l)|]^\top$ and $|\veci(\vecu,\vecs_l)| = [|i_1(\vecu,\vecs_l)|, \dots, |i_L(\vecu,\vecs_l)|]^\top$. In the remainder of this paper, we  proceed under the following practical assumption. 

\begin{assumption}[\textit{Neighborhood of nominal voltage}]
\label{as:steadyStateMap}
\normalfont
The functions $(\vecu, \vecs_l) \mapsto |\vecv(\vecu, \vecs_l)|$ and $(\vecu, \vecs_l) \mapsto |\veci(\vecu, \vecs_l)|$ are unique and continuously differentiable 
in an open neighborhood of the nominal voltage
profile. Additionally, their Jacobian matrices $\boldsymbol{J}_v(\vecu, \vecs_l):=\frac{\partial |\vecv(\vecu, \vecs_l)|}{\partial \vecu}$ and $\boldsymbol{J}_i(\vecu, \vecs_l):=\frac{\partial |\veci(\vecu, \vecs_l)|}{\partial \vecu}$ are locally Lipschitz continuous over that neighborhood.  \hfill $\Box$
\end{assumption}

This assumption is supported by the  findings in, e.g.,~\cite{bolognani2015existence,wang2017existence,bernstein2018load}. This assumption will be utilized only in the analysis of the algorithms; it will not play a role in the algorithmic design and  implementations of the proposed methods.

\subsection{AC OPF Formulation}
\label{subsec:opf-formulations}

Several formulations for the AC OPF at the distribution level have been proposed in the literature; see, for example, the survey~\cite{molzahn2017survey} and the representative works~\cite{bolognani2014distributed,gan2014convex,dall2016optimal,picallo2022adaptive}. In this section, we structure our presentation around an AC OPF formulation that includes constraints on node voltages, line currents, and
 operating ranges of  DERs.

Recall that $\vecu = [p_1, \dots, p_G, q_1, \dots, q_G]^\top \in \mathbb{R}^{2G}$ represents the vector of DER active and reactive power injections, and recall that  \( \mathcal{M} \subseteq \mathcal{N} \) is a set of nodes where  voltages are monitored and controlled. In particular, for the latter, let lower and upper bounds on the voltage magnitudes be denoted as $\underline{V}$ and $\overline{V}$, respectively. Additionally, let $\overline{I}$ be an ampacity limit for the lines  that are monitored.  We then consider the following problem formulation to compute the DERs' power setpoints: 
\begin{align}\label{eq:AC-OPF}
    \textsf{U}^*(\vecs_l,\vectheta) := \arg \min_{\vecu \in \mathcal{C}} \quad & \phi_v(|\vecv(\vecu; \vecs_l)|) + \phi_p(\vecu) \notag
    \\
    \text{s.t.} \quad & \underline{V} \leq |v_j(\vecu; \vecs_l)| \leq \overline{V}, ~~\forall j \in \mathcal{M}\notag
    \\
    & \|i_j(\vecu; \vecs_l)| \leq \overline{I}, ~~~~~~\forall j \in \mathcal{L}\\
    & (p_i, q_i) \in \mathcal{C}_i(\vectheta_{u,i}), ~~ \forall i = 1, \dots, G, \notag
\end{align}
where $\phi_v: \mathbb{R}^{M} \rightarrow \mathbb{R}$,  $\phi_p: \mathbb{R}^{2G} \rightarrow \mathbb{R}$, the set $\mathcal{C}_i(\vectheta_{u,i})$ encodes  constraints for the $i$th DER, and the inequalities in the voltage and current constraints are taken entry-wise. 
The cost $\phi_p(\vecu)$ models DER-related considerations such as minimizing active-power curtailment of renewable resources, limiting excessive reactive power provision to reduce inverter losses, or incorporating economic incentives and fairness criteria among DERs. The voltage-related cost $\phi_v(|\vecv|)$ serves as a performance objective and complements the hard voltage constraints. While the constraints enforce admissible voltage limits, $\phi_v(|\vecv|)$ allows the operator to penalize deviations from a desired voltage profile (e.g., nominal voltages) or promote flatter voltage profiles across the network.

We allow a parametric representation of the set through parameters $\vectheta_{u,i}$; to this end, we assume the set $\mathcal{C}_i(\vectheta_{u,i})$ can be expressed as
\begin{align}
    \mathcal{C}_i(\vectheta_{u,i}) = \{(p_i,q_i) \in \mathbb{R}^2: \ell_i( p_i,q_i, \vectheta_{u,i}) \leq \mathbf{0}_{n_{c_i}}\}
\end{align}
where $\ell_i$ is a vector-valued function modeling power limits, and the inequality is taken entry-wise.  The function $\ell_i$ is assumed to be differentiable. For example,
if the $i$th DER is an inverter-interfaced controllable renewable source, then $\ell_i(p_i,q_i,\vectheta_{u,i}) = [p_i^2+q_i^2 - s_{n,i}^2, p_i - p_{\text{max},i}, -p_i]^\top$, where $\vectheta_{u,i} = (p_{\text{max},i},s_{n,i})$ with $s_{n,i}$ and $p_{\text{max},i}$ the inverter rated size and the maximum available active power, respectively.
The overall  set of inputs that parametrize the problem~\eqref{eq:AC-OPF} is denoted as $\vectheta := (\vectheta_{u,i}, \dots, \vectheta_{u,G}, \underline{V}, \bar{V}, \bar{I})$; these inputs are in addition to $\vecs_l$. We will use the notation $\mathcal{C} = \mathcal{C}_1  \times \ldots \times \mathcal{C}_G$ and we define the set $\mathcal{S}(\vecs_l) := \mathcal{S}_v(\vecs_l) \cap \mathcal{S}_i(\vecs_l)$, where:
\begin{align*}
    \mathcal{S}_v(\vecs_l) & := \{\vecu \in \mathcal{C}: \underline{V} \leq |v_j(\vecu; \vecs_l)| \leq \overline{V}, j\in\mathcal{M}\} \\
    \mathcal{S}_i(\vecs_l) & := \{\vecu \in \mathcal{C}: |i_j(\vecu; \vecs_l)| \leq \overline{I}, j\in\mathcal{L}\} \, .
\end{align*}
The feasible set of~\eqref{eq:AC-OPF} is $\mathcal{S}_v(\vecs_l) \cap \mathcal{S}_i(\vecs_l)$. In the following, for notational simplicity, we drop the dependence on $\vecs_l$.

It is well known that the AC OPF is nonconvex and may admit multiple globally optimal and locally optimal solutions. Accordingly, the function $(\vecs_l,\vectheta) \mapsto \textsf{U}^*(\vecs_l,\vectheta)$ that \emph{maps parameters of the problem into globally optimal solutions to the AC OPF} is a \emph{set-valued function}. 
Since identifying a solution $\vecu^* \in \textsf{U}^*(\vecs_l,\vectheta)$ is in general difficult, we consider a set of local minimizers and isolated KKT points for~\eqref{eq:AC-OPF}; we denote such set as $\textsf{U}^{\textsf{lm}}(\vecs_l,\vectheta)$  (although we note that some local minimizers can also be global minimizers). In the following, we explain our approach to identify  local minimizers in $\textsf{U}^{\textsf{lm}}(\vecs_l,\vectheta)$. 

\subsection{Proposed Mathematical Framework and Implementations}

Our proposed technical approach is grounded on a mathematical model of the form: 
\begin{subequations}
\label{eq:main_method}
\begin{align}
    \dot{\vecu} & = \eta F(\vecu, \vecxi,\vectheta), \label{eq:main_method_flow} \\
    \underbrace{\begin{bmatrix}
       \vecV \\
       \vecI
    \end{bmatrix}}_{:=\vecxi}
     & = \underbrace{\begin{bmatrix}
       |\vecv(\vecu; \vecs_l)| \\
       |\veci(\vecu; \vecs_l)|
    \end{bmatrix}}_{:= H(\vecu; \vecs_l)}  
    + \underbrace{\begin{bmatrix}
       \vecn_v \\
       \vecn_i
    \end{bmatrix}}_{:=\vecn} \label{eq:main_method_pf}
\end{align}
\end{subequations}
where: \emph{(a)} $F(\vecu, \vecxi,\vectheta)$ is a given algorithmic map, utilized to seek local minimizers in $\textsf{U}^{\textsf{lm}}(\vecs_l,\vectheta)$; this map updates $\vecu$ based on voltages $\vecv$, currents $\veci$, and the problem parameters $\vectheta = (\vectheta_{u,i}, \dots, \vectheta_{u,G}, \underline{V}, \bar{V}, \bar{I})$. \emph{(b)}~$H(\vecu; \vecs_l)$ in~\eqref{eq:main_method_pf} represents a  power flow solution map; in particular, given $\vecu, \vecs_l$, one solves for~\eqref{eqn:power-flow} to obtain voltages and currents (and then computes their absolute values). In~\eqref{eq:main_method_pf}, $\vecn$ represents an error or a perturbation in the computation of $\vecxi$. 

\begin{figure*}
    \centering
    \includegraphics[height=5.0cm]{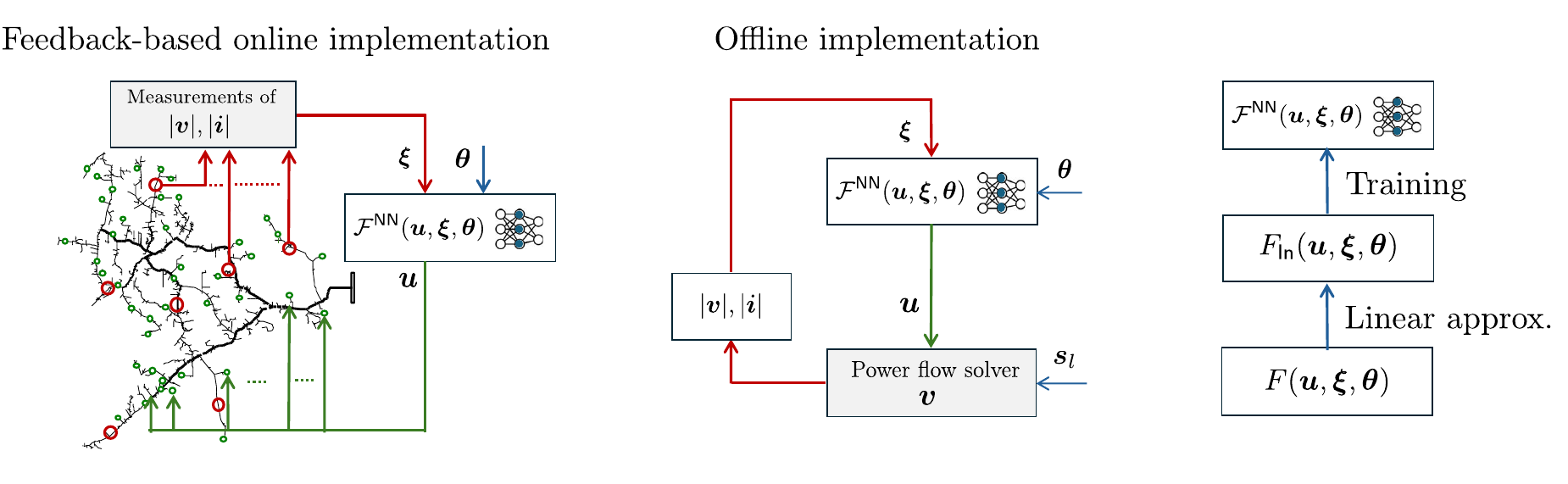}
    \vspace{-.2cm}
    \caption{\emph{(Left)} Feedback-based online implementation, leveraging measurements of voltages and currents from the power distribution network. \emph{(Center)} Offline implementation, where in the iterative process voltages and currents are computed using a power-flow solver. In both cases, $\xi$ collects the voltages and line-currents $(|\vecv|,|\veci|)$. \emph{(Right)} Design process, where the linear approximation of the power flow equations is utilized to lower the computational complexity of the training.}
    \label{fig:proposed_arch}
    \vspace{-.5cm}
\end{figure*}

We design $F(\vecu, \vecxi,\vectheta)$ based on 
the safe gradient flow (SGF)~\cite{allibhoy2023control,colot2024optimal}. In particular, $F(\vecu, \vecxi,\vectheta)$ is given by:
\begin{align}
    & \dot{\vecu} = \eta F(\vecu, \vecxi,\vectheta) \label{eq:controller} \\
    & F(\vecu, \vecxi,\vectheta) := \nonumber\\
    & \arg \min_{\vecz \in \mathbb{R}^{2G}}  \|\vecz  + \nabla \phi_p(\vecu) + \boldsymbol{J}_v(\vecu;\vecs_l)^\top\nabla \phi_v(\vecV)\|^2 \nonumber  \\
    & \hspace{1cm}  \textrm{s.t.}  -\boldsymbol{J}_v(\vecu;\vecs_l)^\top \vecz  \leq -\beta \left(\mathbf{1} \underline{V}- \vecV\right) \label{eq:controller_ideal} \\
    & \hspace{1.8cm} \boldsymbol{J}_v(\vecu;\vecs_l)^\top \vecz  \leq -\beta \left(\vecV -\bar{V} \mathbf{1}\right) \nonumber \\  & \hspace{1.8cm} \boldsymbol{J}_i(\vecu;\vecs_l)^\top \vecz \leq -\beta \left(\vecI -\bar{I} \mathbf{1}\right) \nonumber \\
    &  \hspace{1.7cm} \boldsymbol{J}_{\ell_i}(p_i,q_i)^\top \vecz  \leq - \beta \ell_i(p_i,q_i),  \qquad   \forall i \in \mathcal{G} \nonumber 
\end{align}
where $\boldsymbol{J}_{\ell_i}(p_i,q_i)$ is the Jacobian of $(p_i,q_i) \mapsto \ell_i(p_i,q_i)$,  $\beta > 0$ is a design parameter, and $\eta > 0$ is the controller gain and is a design parameter (inequalities are taken entry-wise). As discussed in~\cite{allibhoy2023control}, the controller in~\eqref{eq:controller} serves as an approximation of the projected gradient flow. This approximation, which leverages control barrier function (CBF) tools, ensures that the feasible set of~\eqref{eq:AC-OPF} is forward invariant. This invariance property is a key motivation behind initiating our design from~\eqref{eq:controller_ideal}, and it is supported by our recent work in~\cite{colot2024optimal}, where the function $F(\vecu, \vecxi, {\vectheta})$ was specifically designed to address an AC OPF problem with voltage constraints.

On the other hand, the update in~\eqref{eq:main_method_pf} lends itself to two distinct practical implementations as outlined below:

\noindent $\triangle$ \emph{Online feedback-based implementation:} Once the update~\eqref{eq:main_method_flow} is performed, the power setpoints $\vecu$ are transmitted to (and  implemented by) the DERs; then, the system operator collects measurements of actual voltages and currents from the system, or leverages pseudo-measurements. This implementation is aligned with existing works on feedback-based optimization~\cite{bolognani2014distributed,gan2016online,dall2016optimal,bernstein2019real,picallo2022adaptive,colot2024optimal}, and it is illustrated in Figure~\ref{fig:proposed_arch}.  

\noindent $\triangle$ \emph{Model-based offline implementation:} In this case,~\eqref{eq:main_method_pf} represents the solution to the AC power flow equations~\eqref{eqn:power-flow} via numerical methods. For example, given $\vecs_l$, voltages $\vecv$ can be found using the fixed-point method~\cite{bernstein2018load,chen2023fixed}. This leads to the offline solution of~\eqref{eq:AC-OPF} illustrated  in Figure~\ref{fig:proposed_arch}.      

When $F(\vecu, \vecxi,\vectheta)$ involves the solution of a quadratic program (QP) as in~\eqref{eq:controller_ideal}, the  process in~\eqref{eq:main_method} can be computationally intensive; the time required to identify a solution may not align with the time scales at which loads $\vecs_l$ and system parameters $\vectheta$ evolve~\cite{taylor2016power}. More broadly, similar arguments would apply to cases where $F(\vecu, \vecxi,\vectheta)$ is designed using different algorithmic approaches involving 
projections onto manifolds~\cite{Hauswirth16} or inversions of (potentially large) matrices as in  Newton-type methods~\cite{baker2020learning}. The idea is then to train a neural network to approximate $F(\vecu, \vecxi,\vectheta)$. Letting $\mathcal{F}^{\textsf{NN}}: \mathbb{R}^{2G} \times \mathbb{R}^{M+L} \times \mathbb{R}^{n_\theta} \rightarrow \mathbb{R}^{2G}$ the neural network map, we consider the following modification of~\eqref{eq:main_method}:
\begin{subequations}
\label{eq:main_NN}
\begin{align}
    \dot{\vecu} & = \eta \mathcal{F}^{\textsf{NN}}(\vecu,\vecxi,\vectheta) , \label{eq:main_NN_flow} \\
    \vecxi & = H(\vecu; \vecs_l) + \vecn \label{eq:main_NN_pf}
\end{align}
\end{subequations}
where we recall that $\vecxi$ represents either a solution to the power flow equations via numerical methods or measurements of voltages and currents. Based on the model~\eqref{eq:main_NN} the problem addressed in the remainder of the paper is as follows.

\begin{problem}
\label{prob:main}
    Train a neural network $\mathcal{F}^{\textsf{NN}}$ to emulate $F(\vecu, \vecxi,\vectheta)$ in~\eqref{eq:controller_ideal} so that the algorithm~\eqref{eq:main_NN}: (a)~converges to a solution $\vecu^* \in \textsf{U}^{\textsf{lm}}(\vectheta)$   of the AC OPF problem \eqref{eq:AC-OPF}; (b)~ensures that voltage, current, and DER  constraints are satisfied at any time during the execution of the algorithm.  \hfill $\Box$
\end{problem}

The term ``at any time'' refers to the fact that the algorithm~\eqref{eq:main_NN} is expected to produce points that are practically feasible for~\eqref{eq:AC-OPF} even if it is terminated before convergence. This is a key for online AC OPF implementations as in Figure~\ref{fig:proposed_arch}(left), and a desirable feature of offline methods.  We provide some remarks to support our design approach.

\subsection{Motivations and Rationale}

An approach different than the one proposed in~\eqref{eq:main_NN} in the context of learning for the AC OPF problems is to train a neural network to directly identify optimal solutions in $\textsf{U}^*(\vecs_l,\vectheta)$, solutions in $\textsf{U}^{\textsf{lm}}(\vecs_l,\vectheta)$, or the set of KKT points $\textsf{U}^{\textsf{kkt}}(\vecs_l,\vectheta)$; see, for example,~\cite{zamzam2020learning,pan2023deepopf,park2023compact,fioretto2020predicting} and \cite{hamann2024foundation}. We provide a comparison next.

\noindent $\triangle$ \emph{Mapping vs set-valued mapping} 
\begin{itemize}
    \item  For any given $\vecu, \vecxi,\vectheta$, $F(\vecu, \vecxi,\vectheta)$ is defined as the \emph{unique} optimal solution of the convex QP~\eqref{eq:controller_ideal}. Moreover, under some mild assumptions, the mapping $F(\vecu, \vecxi,\vectheta)$ is locally Lipschitz (in all its arguments)~\cite{allibhoy2023control,colot2024optimal}. Therefore, in our approach,  the neural network approximates a mapping that is continuous in its arguments. 
    \item On the other hand, $\textsf{U}^*(\vecs_l,\vectheta)$, $\textsf{U}^{\textsf{lm}}(\vecs_l,\vectheta)$, and $\textsf{U}^{\textsf{kkt}}(\vecs_l,\vectheta)$ are in general \emph{sets}; therefore, $(\vecs_l,\vectheta) \mapsto \textsf{U}^*(\vectheta)$ is a set-valued mapping. Such a set-valued mapping cannot be approximated with the mapping of the neural network; see the discussion in, e.g.~\cite{sun2018learning} and~\cite{pan2023deepopf}. One workaround suggested in~\cite{sun2018learning}  is when the number of solutions to the AC OPF is finite and all the solutions can be identified; in this case, letting as an example $\vecu^{\textsf{kkt}}(\vecs_l,\vectheta)$ be a vector collecting all the KKT points for given parameters $(\vecs_l,\vectheta)$, one can use a neural network to approximate $(\vecs_l,\vectheta) \mapsto \vecu^{\textsf{kkt}}(\vecs_l,\vectheta)$. However, the  solutions of the AC OPF cannot be, in general, enumerated.     
\end{itemize}
\noindent $\triangle$ \emph{Number of inputs}
\begin{itemize}
    \item {To approximate $\mathcal{F}^{\textsf{NN}}(\vecu,\vecxi,\vectheta)$, the inputs to the training are: (i) the monitored voltage and/or line-current magnitudes collected in $\vecxi$,  (ii) the current DER setpoints $\vecu$, and (iii) the parameters $\vectheta$ (see Fig.~\ref{fig:proposed_arch}). 
    We recall that $\vecxi$ is obtained from  measurements in the feedback-based online implementation or from numerical power-flow solutions in the offline implementation. The loads $\vecs_\ell$ are not included in the feature vector used for training.}
    \item To approximate $\textsf{U}^*(\vecs_\ell,\vectheta)$, $\textsf{U}^{\textsf{lm}}(\vecs_\ell,\vectheta)$, or $\textsf{U}^{\textsf{kkt}}(\vecs_\ell,\vectheta)$ directly, the inputs are the system-wide loads $\vecs_\ell$. {Thus, when the number of load buses is large (e.g., larger than the number of controlled DERs and the number of monitored voltage/current measurements), the dimension of the corresponding input vector grows with $\dim(\vecs_\ell)$, and can be substantially larger than the input dimension required to learn $F(\vecu,\vecxi,\vectheta)$, which depends on the current setpoints $\vecu$ and on  $\vecxi=\big(|\vecv|,|\veci|\big)$.}
\end{itemize}

\noindent $\triangle$ \emph{Feasibility guarantees}
\begin{itemize}
    \item As shown in Section~\ref{sec:convergence}, our method ensures that iterates $\vecu$ are  practically feasible; i.e., we characterize the worst-case violation of a constraint. With this information, and by tightening the constraints during the training process, one can ensure that our method generates points that are feasible for the AC OPF.  
    \item Existing methods that ``emulate'' solutions in $\textsf{U}^*(\vecs_l,\vectheta)$, $\textsf{U}^{\textsf{lm}}(\vecs_l,\vectheta)$, or $\textsf{U}^{\textsf{kkt}}(\vecs_l,\vectheta)$ do not guarantee feasibility of the generated outputs. Post-processing could adjust the solution to make it feasible, but that may involve  heuristics that do not have feasibility guarantees.   
\end{itemize}


\section{Neural Network-based OPF Pursuit}
\label{sec:NNopf}

In this section, we provide details on the algorithmic design and we discuss its implementation. 

\subsection{Algorithmic Design}

The map $F(\vecu, \vecxi,\vectheta)$ in~\eqref{eq:controller} requires computing the Jacobian matrices of function $H(\vecu, \vecs_l)$. To favor a lower complexity training procedure, we rely on a linear approximation of the power flow equations~\eqref{eqn:power-flow}; several linear approximation approaches can be found in the literature; see for example,~\cite{bolognani2015existence,bernstein2018load,gan2014convex} and references therein. In general,  one can find linear approximations of the form
\begin{align}
|\vecv(\vecu; \vecs_l)| \approx 
& \boldsymbol{\Gamma}_v \vecu  +  \bar{\vecv}(\vecs_l), ~~
|\veci(\vecu; \vecs_l)| \approx \boldsymbol{\Gamma}_i \vecu  +  \bar{\veci}(\vecs_l) , \label{eq:linear-model}
\end{align}
where the matrices $\boldsymbol{\Gamma}_v, \boldsymbol{\Gamma}_i$ and the vectors $\bar{\vecv}(\vecs_l), \bar{\veci}(\vecs_l)$ can be computed using the methods in~\cite{bolognani2015existence,bernstein2018load,gan2014convex}. 
The matrices $\boldsymbol{\Gamma}_v, \boldsymbol{\Gamma}_i$ can be precomputed, as they do not depend on $\vecu$ or $\vecs_l$. Using~\eqref{eq:linear-model}, we can utilize the following approximation of~\eqref{eq:controller_ideal}: 
\begin{align}
    & \dot{\vecu} = \eta F_{\textsf{ln}}(\vecu, \vecxi,\vectheta) \label{eq:controller_appr} \\
    & F_{\textsf{ln}}(\vecu, \vecxi,\vectheta) := \arg \min_{\vecz \in \mathbb{R}^{2G}}  \|\vecz  + \nabla \phi_p(\vecu) + \boldsymbol{\Gamma}_v^\top\nabla \phi_v(\vecV)\|^2 \nonumber  \\
    & \hspace{3cm}  \textrm{s.t.}  -\boldsymbol{\Gamma}_v^\top \vecz  \leq -\beta \left(\mathbf{1} \underline{V}- \vecV\right) \label{eq:controller_approximated} \\
    & \hspace{3.8cm} \boldsymbol{\Gamma}_v^\top \vecz  \leq -\beta \left(\vecV -\bar{V} \mathbf{1}\right) \nonumber \\  & \hspace{3.8cm} \boldsymbol{\Gamma}_i^\top \vecz \leq -\beta \left(\vecI -\bar{I} \mathbf{1}\right) \nonumber \\
    &  \hspace{3.7cm} \boldsymbol{J}_{\ell_i}(\vecu_i)^\top \vecz  \leq - \beta \ell_i(\vecu), ~ i \in \mathcal{G} \nonumber 
\end{align}
where we have replaced the Jacobian matrices of the power flow equations with the ones of the linear approximations. In our forthcoming analysis in Section~\ref{sec:convergence}, we quantify the effect of the linear approximation error in the overall performance. 

Similarly to~\eqref{eq:controller_ideal}, ~\eqref{eq:controller_approximated} is a convex QP with a unique optimal solution. 
Moreover, from~\cite[Theorem 3.6]{liu1995sensitivity}, it follows that $\vecu \mapsto F_{\textsf{ln}}(\vecu, \vecxi,\vectheta)$ is locally Lipschitz over $\mathcal{B}(\vecu,r_1):=\{\vecz : \|\vecz-\vecu\|<r_1 \}$ of $\vecu$, for any $\vecxi$ and $\vectheta$.

The next step, as illustrated in Figure~\ref{fig:proposed_arch}(right), is to consider a neural network $\mathcal{F}^{\textsf{NN}}: \mathbb{R}^{2G} \times \mathbb{R}^{M+L} \times \mathbb{R}^{n_\theta} \rightarrow \mathbb{R}^{2G}$, which will be trained to approximate the mapping $(\vecu, \vecxi,\vectheta) \mapsto F_{\textsf{ln}}(\vecu, \vecxi,\vectheta)$. In particular, 
consider a fully connected feedforward neural network (FNN), defined recursively as:
\begin{align}
\vecy &= \mathcal{F}^{\textsf{NN}}(\vecu, \vecxi,\vectheta) := W^{(H)}\boldsymbol{\varphi}^{(H)} + b^{(H)}, \label{eq:nn-final} \\
\boldsymbol{\varphi}^{(i)} &= \Phi^{(i)}\left(W^{(i-1)}\boldsymbol{\varphi}^{(i-1)} + b^{(i-1)}\right), \quad i = 1, \ldots, H, \nonumber \\
\boldsymbol{\varphi}^{(0)} &= [\vecu, \vecxi,\vectheta], \nonumber
\end{align}
where $H$ is the number of hidden layers, $W^{(i)} \in \mathbb{R}^{n_{i+1} \times n_i}$ and $b^{(i)} \in \mathbb{R}^{n_{i+1}}$ are the weights and biases, and $\Phi^{(i)}$ is a Lipschitz-continuous activation function (e.g., ReLU, leaky ReLU, or sigmoid) 
The network outputs in~\eqref{eq:nn-final} are $\vecy = \dot{\vecu}$. 

For the training procedure, suppose that $N_{\text{train}}$ training points are available, and they are taken from a compact set $(\vecu, \vecxi, \vectheta) \in \mathcal{C}_{\text{train}} \times \mathcal{E}_{\text{train}} \times \Theta_{\text{train}}$, where $\mathcal{C}_{\text{train}}$ is a superset of the feasible region of~\eqref{eq:AC-OPF}, $\mathcal{E}_{\text{train}}$ is an inflation of the set of operational voltages, and $\Theta_{\text{train}}$ is formed based on inverters' operating conditions. Thus, each training point is given by the input $(\vecu^{(k)}, \vecxi^{(k)}, \vectheta^{(k)})$ and the corresponding output $\vecy^{(k)} = F_{\textsf{ln}}(\vecu^{(k)}, \vecxi^{(k)}, \vectheta^{(k)})$, for $k = 1, \ldots, N_{\text{train}}$. Then, we consider minimizing the following loss function: 
\begin{align}
\hspace{-.2cm} \mathcal{L}(\boldsymbol{W},\boldsymbol{b}) := \frac{1}{N_{\text{train}} } \sum_{n=1}^{N_{\text{train}}}
\left\|\vecy^{(k)} - \mathcal{F}^{\mathsf{NN}}(\vecu^{(k)}, \vecxi^{(k)}, \vectheta^{(k)}) \right\|_2^2
\end{align}
where the dependence of $\mathcal{F}^{\mathsf{NN}}$ on $\boldsymbol{W},\boldsymbol{b}$ is dropped for notational convenience. The training routine is presented in Algorithm~\ref{alg:nn-sgf-offline}. 

\begin{algorithm}[H]
\caption{Offline training}
\label{alg:nn-sgf-offline}
\begin{algorithmic}[1]
\STATE \textbf{Generate or collect training points}
\STATE ~~~ For each time instant or episode \( \{t_k\}_{k=1}^{N_{\text{train}}}\):
\STATE ~~~~~ Obtain $\vecv^{(k)}$, $\veci^{(k)}$, and $\vecu^{(k)}$
\STATE ~~~~~ Obtain parameters: \( \vectheta^{(k)} \)
\STATE ~~~~~ Compute: $\vecy^{(k)} = F_{\text{ln}}(\vecu^{(k)}, \vecxi^{(k)}, \vectheta^{(k)})$
\STATE \textbf{Train neural network}
\STATE ~~~ Solve $\min_{\boldsymbol{W},\boldsymbol{b}} \mathcal{L}(\boldsymbol{W},\boldsymbol{b})$
\end{algorithmic}
\end{algorithm}

We note that the training dataset can be generated offline by repeating {the data-generation loop in Algorithm~\ref{alg:nn-sgf-offline} (steps 1--5)} for a given set of values for voltages, currents, and DERs' powers, or it can be formed online by collecting measurements from the distribution grid. The trained  FNN $\mathcal{F}^{\mathsf{NN}}(\vecu, \vecxi,\vectheta)$ is then used in~\eqref{eq:main_NN} to solve the AC OPF online (cf. Figure~\ref{fig:proposed_arch}(left) or offline (cf. Figure~\ref{fig:proposed_arch}(center)). 


\subsection{Online and Offline Implementations}

In this section, we provide more details on the online and offline implementations of our proposed method. The feedback-based online implementation is illustrated in Figure~\ref{fig:proposed_arch}(left); here, the parameters $\vectheta(t)$ are time-varying since they include the power available from renewable-based DERs, which may change with evolving ambient conditions~\cite{dall2016optimal}. The overall algorithm is tabulated as Algorithm~\ref{alg:nn-online}. Similar to existing feedback-based algorithms, Algorithm~\ref{alg:nn-online} does not require any information about the loads $\vecs_l$. In this implementation, the error term $\vecn$ in~\eqref{eq:main_method_flow} 
 and~\eqref{eq:main_NN_flow} represents errors in the measurements of voltages and currents, or in the computation of pseudo-measurements; these errors are small or even negligible~\cite{angioni2015impact}. 

On the other hand, the offline implementation of  Figure~\ref{fig:proposed_arch}(center) is tabulated as Algorithm~\ref{alg:nn-offline}. For this offline implementation, one needs information about the loads $\vecs_l$. {In Algorithm~\ref{alg:nn-offline} (step~5),} solutions to the power flow (PF) equations can be identified using, for example, sweeping methods~\cite{Kerstingbook} or fixed-point methods~\cite{bernstein2018load}. In this implementation, the error $\vecn$ in~\eqref{eq:main_method_flow} 
 and~\eqref{eq:main_NN_flow} represents the numerical accuracy of the PF method. The algorithm is executed until convergence, or for a pre-scribed amount of time~$t_d$. 

\begin{algorithm}[t!]
\caption{Online feedback-based implementation}
\label{alg:nn-online}
\begin{algorithmic}[1]
\STATE \textbf{Initialization}
\STATE ~~~ Load pretrained model \( \mathcal{F}^{\mathsf{NN}}\), pick $\eta > 0$.
\STATE \textbf{Real-Time operation \( t \geq 0 \):}
\STATE ~~~ Measure DERs' setpoints \( \vecu(t)\)
\STATE ~~~ Measure $|\vecv(t)|$ and $|\veci(t)|$ from selected locations
\STATE ~~~ Obtain parameters $\vectheta(t)$
\STATE ~~~ Perform update:
$\dot{\vecu}(t) = \eta \mathcal{F}^{\mathsf{NN}} \left( \vecu(t), \vecxi(t), \vectheta(t) \right)$
\STATE ~~~ Send \( \vecu(t) \) to DERs  and go to 4.
\end{algorithmic}
\end{algorithm}

\section{Experimental Results in a Distribution Feeder}
\label{sec:results}
\begin{figure}[t]
\centering 
    \includegraphics[height=4.6cm]{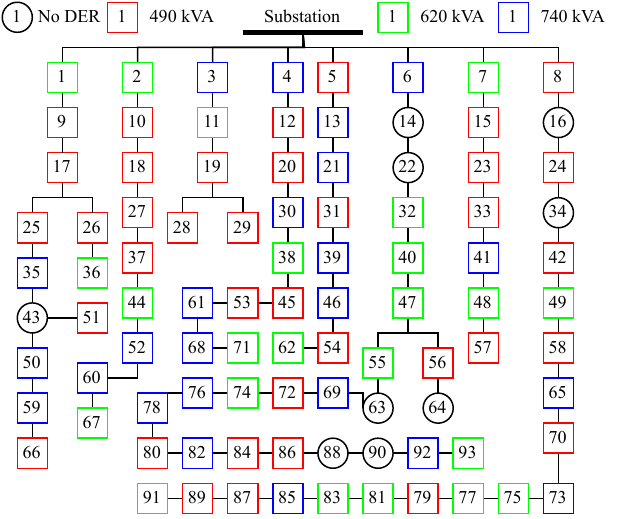}
\caption{Distribution network used in the simulations~\cite{SimBench}.}
\label{fig:network}
\vspace{-.2cm}
\end{figure}

We test the proposed method illustrated in Figure~\ref{fig:proposed_arch}(left) -- which in this section we refer to as neural network-based safe gradient flow (NN-SGF) in short - on a voltage regulation problem. We consider the medium voltage network (20~kV) shown in Figure~\ref{fig:network} (see~\cite{SimBench}). {We define all per-unit quantities using $V_{\text{base}} = 20~\text{kV}$ and $S_{\text{base}} = 10~\text{MVA}$, yielding $Z_{\text{base}} = V_{\text{base}}^2/S_{\text{base}} = 40~\Omega$.}
The network contains photovoltaic~(PV) inverters at selected buses capable of adjusting both active and reactive power. Each inverter \( i \in \mathcal{G} \) injects power \( \vecu_i = (p_i, q_i) \in \mathbb{R}^2 \) within a feasible set:
\begin{equation}
\label{eq:Ui}
\mathcal{C}_i(\vectheta_{u,i}) = \left\{ (p_i, q_i) \in \mathbb{R}^2 \;\middle|\;
\begin{bmatrix}
p_i^2 + q_i^2 - s_{n,i}^2 \\
p_i - p_{\max,i} \\
- p_i \\
-0.44\,s_{n,i} - q_i \\
q_i - 0.44\,s_{n,i}
\end{bmatrix}
\leq \mathbf{0}
\right\},
\end{equation}
where \(s_{n,i}\) represents the inverter's nominal apparent power rating, randomly selected from the set \{490, 620, 740\}~kVA {(i.e., \(s_{n,i}/S_{\text{base}} \in \{0.049,\,0.062,\,0.074\}\) p.u.)} to capture the range of deployment scales observed in practice. The upper bound \(p_{\max,i}\) denotes the maximum available active power at time \(t\). Together, the pair \( \vectheta_{u,i} = (p_{\max,i}, s_{n,i}) \) defines the parameterization of the  set \( \mathcal{C}_i(\vectheta_{u,i}) \).
This limit is consistent with practical deployment settings and in accordance with IEEE Std 1547-2018. We also assume that \(p_{\max,i}\) is known at the DER level (via maximum power point tracking).  The cost function in the AC OPF  is defined as $\phi_p(\vecu) = \sum_{i \in \mathcal{G}} 
c_p \left( \frac{s_{n,i} - p_i}{s_{n,i}} \right)^2 + 
c_q \left( \frac{q_i}{s_{n,i}} \right)^2$, with \(c_p = 3\) and \(c_q = 1\). Here, the first term aims to minimize active power curtailment and inverter losses, whereas the latter discourages reliance on reactive power, as it leads to increased current levels and energy losses.

The voltage magnitudes at monitored buses are denoted by \(\vecV\) (cf.~\eqref{eq:main_method_pf}) and are obtained from the \texttt{pandapower} power flow solver\footnote{See \url{https://www.pandapower.org}.}. The aggregated non-controllable loads and maximum available active power for PV plants is from the Open Power System 
Data\footnote{\url{https://data.open-power-system-data.org/household_data/2020-04-15}}; the data has a granularity of 10 seconds, and the values have been modified to match the initial loads and PV plants nominal values present in the network. The reactive power demand is set such that the power factor is 0.9 (lagging). The voltage service limits $\bar{V}$ and $\underline{V}$ are set to 1.05 and 0.95 p.u., respectively. 

\begin{algorithm}[t!]
\caption{Offline implementation}
\label{alg:nn-offline}
\begin{algorithmic}[1]
\STATE \textbf{Initialization}
\STATE ~~~ Load pretrained model \( \mathcal{F}^{\mathsf{NN}}\), pick $\eta > 0$.
\STATE ~~~ Load parameters $\vectheta$, load $\vecs_l$, set $\vecu(0)$.
\STATE \textbf{Perform until convergence or until $t_d$:}
\STATE ~~~ Given $\vecs_l, \vecu(\tau)$,  solve PF equations to get $\vecv(\tau)$
\STATE ~~~ Compute $|\vecv(\tau)|$ and $|\veci(\tau)|$ from selected locations
\STATE ~~~ Perform update:
$\dot{\vecu}(\tau) = \eta \mathcal{F}^{\mathsf{NN}} \left( \vecu(\tau), \vecxi(\tau), \vectheta \right)
$
\STATE ~~~ Go to 5.
\end{algorithmic}
\end{algorithm}

\begin{figure}[t]
    \centering
    \includegraphics[height=4.6cm]{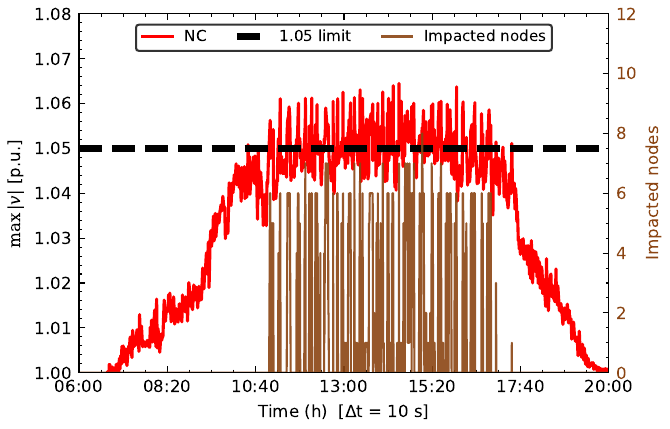}
    \caption{Overvoltage events under the considered simulation setup when no control actions are implemented. {Evaluated over the monitored buses $\mathcal{M}=\{93,54,85,84,48,39,91,71,62,67\}$.}}
    \label{fig:voltage_nc}
\vspace{-.4cm}
\end{figure}

With the considered data and simulation setup, we obtain the voltage profiles illustrated in Figure~\ref{fig:voltage_nc} for the case of \emph{no control}; this is a case where a protection scheme of the PV plants disconnects the inverters if the voltage level is too high. The disconnection scheme is inspired from the CENELEC EN50549-2 standard; the PV plant changes status from running to disconnected if: (i) the voltage at the point of connection goes above 1.06 pu, (ii) the root mean square value of the voltages measured at the point of connection for the past 10 minutes goes above 1.05 pu (the voltages are measured every 10 seconds). The switching from disconnected to connected occurs randomly
in the interval [1min, 10min]. 

As shown in Figure~\ref{fig:voltage_nc}, the proposed method is tested against a challenging voltage regulation problem. We compare the solutions obtained with the following strategies: 

\noindent $\triangle$ \emph{(s1)} A solution of the AC OPF every 10 seconds to match the granularity of the Open Power System Data. Here, we use the nonlinear branch flow model \cite{BFM} and the solver \texttt{IPOPT}. We refer to this case as batch optimization (BO), {meaning that at every time we solve the nonlinear AC OPF problem  to convergence.}

\noindent $\triangle$ \emph{(s2)} Our solution strategy in~\eqref{eq:main_NN} is deployed in the online feedback-based configuration in Figure~\ref{fig:proposed_arch} (left). Here, we run \emph{one iteration} of NN--SGF each time a measurement is collected (as in standard feedback optimization methods). {Here, ``one iteration'' means a single forward evaluation of $\mathcal{F}^{\textsf{NN}}\!\big(\boldsymbol{u}(t_k),\boldsymbol{\xi}(t_k),\boldsymbol{\theta}(t_k)\big)$ using the latest measurement at time $t_k$; the setpoints are then held  until the next measurement arrives.}

\noindent $\triangle$ \emph{(s3)} A strategy similar to, e.g.,~\cite{zamzam2020learning,nellikkath2022physics,pan2022deepopf} where a neural network is used to approximate $\textsf{U}^{\textsf{lm}}(\vecs_l,\vectheta)$; i.e., to emulate solutions of the BO directly.  We refer to this as NN-BO.

\subsection{Dataset Generation and Training}

We simulate \( N_{\text{training}} = 6{,}000 \) independent operating conditions, each defined by a distinct grid configuration with varying load profiles, DER capacities, and voltage regulation constraints. Each condition is associated with a randomly sampled time instant \( t_n \sim \mathcal{U}(t_{\min}, t_{\max}) \),
%
%
where \( t_{\min} = \text{06{:}00} \), \( t_{\max} = \text{20{:}00} \). For each sampled time \( t_n \), we simulate power flow using \texttt{pandapower} to obtain voltages \( \vecV(t_n) = \{ V_j(t_n) \}_{j \in \mathcal{M}} \), along with DER setpoints \( \{ p_i(t_n), q_i(t_n) \}_{i \in \mathcal{G}} \), constraint parameters \( \{ p_{\max,i}(t_n), s_{n,i} \}_{i \in \mathcal{G}} \), and voltage bounds \( \{ \underline{V}_j, \overline{V}_j \}_{j \in \mathcal{M}} \). For each operating condition, we run \( N_{\text{iter}} = 10 \) iterations of \( F_{\textsf{ln}}(\vecu, \vecxi, \vectheta) \) using a forward Euler discretization, with \( \eta = 0.2 \). At each iteration \( k \in \{1, \ldots, N_{\text{iter}}\} \), we record the DER setpoints \( \vecu^{(n,k)} = \{ p_i^{(k)}(t_n), q_i^{(k)}(t_n) \}_{i \in \mathcal{G}} \), the voltage measurements \( \vecV^{(n,k)} = \{ V_j^{(k)}(t_n) \}_{j \in \mathcal{M}} \), and the constraint parameters \( p_{\max,i}(t_n) \), \( s_{n,i} \), \( \underline{V}_j \), and \( \overline{V}_j \). We also extract the current magnitudes \( \vecI^{(n,k)} \in \mathbb{R}^L \) at monitored lines. These are used to construct the state vector \( \vecxi^{(n,k)} = [ (\vecV^{(n,k)})^\top,\, (\vecI^{(n,k)})^\top ]^\top \in \mathbb{R}^{M+L} \), consistent with the online and offline implementation setup. The full constraint parameter vector is denoted \( \vectheta^{(n)} = (\vectheta_{u,1}^{(n)}, \ldots, \vectheta_{u,G}^{(n)}, \underline{V}, \overline{V}, \overline{I}) \), and is treated as fixed across SGF iterations at time \( t_n \). The control update is computed as \( \vecy^{(n,k)} = F_{\textsf{ln}}(\vecu^{(n,k)}, \vecxi^{(n,k)}, \vectheta^{(n)}) \in \mathbb{R}^{2G} \) and used as the training label. The corresponding input vector \( \mathbf{x}^{(n,k)} \in \mathbb{R}^{3G + M} \) is constructed by concatenating the normalized active power deviations \( \left\{ (p_i^{(k)}(t_n) - p_{\max,i}(t_n))/s_{n,i}, i \in \mathcal{G} \right\} \), normalized reactive powers \( \left\{q_i^{(k)}(t_n)/s_{n,i}, i \in \mathcal{G} \right\} \), normalized voltage magnitudes \( \left\{ (V_j^{(k)}(t_n) - \underline{V}_j)/(\overline{V}_j - \underline{V}_j), j \in \mathcal{M} \right\} \), and the active DER limits \( \left\{ p_{\max,i}(t_n), i \in \mathcal{G} \right\} \). This normalization ensures all features lie in comparable ranges and are scaled relative to their physical limits (e.g., \( s_{n,i} \) and voltage bounds), improving numerical stability. This process results in a dataset of \( N_{\text{train}} = N_{\text{training}} \cdot N_{\text{iter}} = 60{,}000 \) input-output pairs \( \mathcal{D}_{\text{train}} = \left\{ \left( \mathbf{x}^{(n,k)}, \mathbf{y}^{(n,k)} \right) \right\}_{n,k} \subset \mathbb{R}^{3G+M} \times \mathbb{R}^{2G} \). For evaluation, we construct a disjoint test set \( \mathcal{D}_{\text{test}} = \left\{ \left( \mathbf{x}^{(n,k)}, \mathbf{y}^{(n,k)} \right) \right\}_{n,k} \subset \mathbb{R}^{3G+M} \times \mathbb{R}^{2G} \) using \( N_{\text{testing}} = 1000 \) randomly sampled times \( t_n \sim \mathcal{U}(\text{06{:}00},\;\text{20{:}00}) \), with the same SGF update procedure repeated for each test case. We ensure that \( \mathcal{D}_{\text{train}} \cap \mathcal{D}_{\text{test}} = \emptyset \).

\emph{Learning Model, Training, and Evaluation}. We train a fully connected FNN \( \mathcal{F}^{\textsf{NN}}: \mathbb{R}^{2G} \times \mathbb{R}^{M+L} \times \mathbb{R}^{n_\theta} \rightarrow \mathbb{R}^{2G} \); the FNN has an architecture of the form \( [3G + M, h, h, h, 2G] \), with hidden width set to \( h = \alpha(3G + M) \); {here, $\alpha$ scales the hidden width and thus directly sets the model capacity.} We use \( \alpha = 2 \), yielding \( h = 524 \) for \( G = 84 \) and \( M = 10 \).
The network is implemented in \texttt{PyTorch} and trained offline using the Adam optimizer\footnote{See: \url{https://docs.pytorch.org/docs/stable/generated/torch.optim.Adam.html}} with learning rate 0.001, batch size 256, dropout 0.2, and up to 500 epochs with early stopping based on validation loss from a 10\% held-out subset. {All experiments  were run on a laptop with  a 16-core Intel Core Ultra 7 255H CPU and 32 GB of RAM. Under this setup, training NN--SGF takes approximately 75 minutes.} The loss function is the mean squared error (MSE) \( \mathcal{L}(\theta) = \frac{1}{|\mathcal{D}_{\text{train}}|} \sum_{n,k} \left\| \mathcal{F}^{\textsf{NN}}(\vecu^{(n,k)}, \vecxi^{(n,k)}, \vectheta^{(n)}; \theta) - \vecy^{(n,k)} \right\|_2^2 \) over the dataset \( \mathcal{D}_{\text{train}} = \left\{ (\vecu^{(n,k)}, \vecxi^{(n,k)}, \vectheta^{(n)}, \vecy^{(n,k)}) \right\}_{n,k} \). Performance on a disjoint test set \( \mathcal{D}_{\text{test}} \) is evaluated using the MSE \( \ell_2 \)-norm prediction error, defined as \( \varepsilon^{\textsf{NN}} = \frac{1}{|\mathcal{D}_{\text{test}}|} \sum_{n,k} \left\| \mathcal{F}^{\textsf{NN}}(\vecu^{(n,k)}, \vecxi^{(n,k)}, \vectheta^{(n)}) - \vecy^{(n,k)} \right\|_2^2 \) ,yielding \( \varepsilon^{\textsf{NN}} = 1.7 \times 10^{-6} \) and \( \mathrm{RMSE} = 0.0013 \). During online deployment, the FNN runs in inference mode at each control step \( t \), with sampling interval \( \Delta t = 10 \) seconds (to match the variability of the load and PV data). Given current DER setpoints \( \vecu(t) \), measurements \( \vecxi(t) \), and parameters \( \vectheta(t) \), the update is approximated as \( \vecu(t + \Delta t) = \vecu(t) + \eta \Delta t \cdot \mathcal{F}^{\textsf{NN}}(\vecu(t), \vecxi(t), \vectheta(t)) \), {where $\eta$ is the controller gain of the safe gradient flow controller (higher values lead to more aggressive updates)}, and it is set to \( \eta = 0.02 \); setpoints are restricted to the set \( \mathcal{C} \), via a projection, if not feasible to reflect hardware constraints. Given the setpoints, updated voltage magnitudes are computed via AC power flow using \texttt{pandapower}, yielding the new vector \( \vecxi(t + \Delta t) \) for the next control step.

For the method \emph{(s3)} used for comparison, the training of an FNN $\mathcal{F}_{\textsf{batch}}^{\textsf{NN}} : (\vecs_l, \vectheta) \mapsto \vecu^\ast$ to emulate solutions to the BO was computationally heavier as we had to increase the size of the training set to obtain acceptable performance. {In fact, by using 60,000 training points (as for the NN-SGF), the performance of the NN-BO was poor; we then resorted to 20,000 additional training points 
for the  NN--BO to obtain a smaller test error. This increased number of  data points reflects the higher sample complexity of directly learning the nonconvex AC--OPF solution map.} The training set consists of inputs \( \vecs_l^{(k)}, \vectheta^{(k)} \), with corresponding outputs \( \vecu^{\ast(k)} \) obtained by solving the AC OPF using IPOPT. To generate the dataset, we sample \( N_{\text{cases}} = 8{,}000 \) time instants \( t_n \sim \mathcal{U}(\text{06:00}, \text{20:00}) \), record the uncontrollable loads \( \vecs_l(t_n) := (\vecp_l(t_n), \vecq_l(t_n)) \) and inverter limits \( \vecp_{\max}(t_n) \), and apply perturbations \( \vecs_l^{(n,k)} := (1 + \epsilon_k)\vecs_l(t_n) \) with \( \epsilon_k \sim \mathcal{U}(-0.05, 0.05) \) for \( k = 1,\dots,10 \), introducing up to \( \pm5\% \) variability to enable localized sampling of the solution space without violating OPF feasibility at \( t_n \). Solving the AC OPF for each perturbed point yields the target \( \vecu^{\ast(n,k)} \in \textsf{U}^{\textsf{lm}}(\vecs_l^{(n,k)}, \vectheta)) \), giving \( N_{\text{train}} = 80{,}000 \) training pairs. The loss minimized is $\mathcal{L}(\theta) = \frac{1}{N_{\text{train}}} \sum_{i=1}^{N_{\text{train}}} \left\| \mathcal{F}_{\textsf{batch}}^{\textsf{NN}}(\vecs_l^{(i)}, \vectheta; \theta) - \vecu^{\ast(i)} \right\|_2^2$. {Using the same laptop with a 16-core
Intel Core Ultra 7 255H CPU and 32 GB of RAM, training NN--BO takes approximately $ 110$ minutes.}

\begin{figure*}[t]
\centering
\begin{subfigure}[t]{0.325\textwidth}
    \centering
    \includegraphics[width=\linewidth, height=4.6cm]{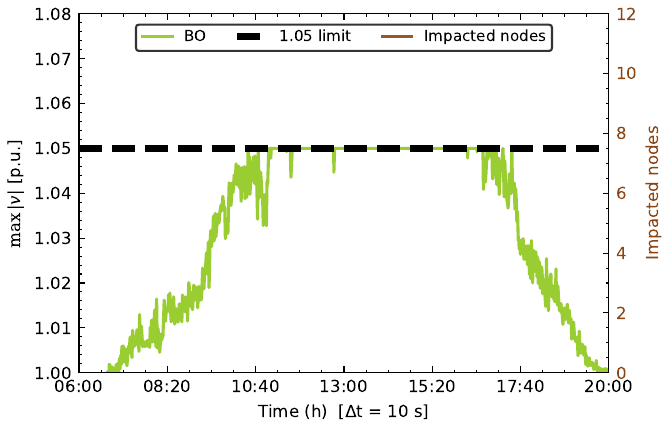}
    \caption{AC OPF solution via IPOPT.}
\end{subfigure}
\hfill
\begin{subfigure}[t]{0.325\textwidth}
    \centering
    \includegraphics[width=\linewidth, height=4.6cm]{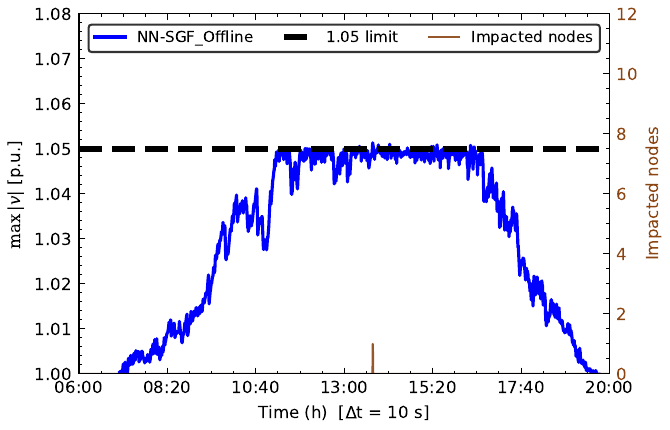}
    \caption{Offline NN-SGF~\eqref{eq:main_NN} in Algorithm~\ref{alg:nn-offline}.}
\end{subfigure}
\hfill
\begin{subfigure}[t]{0.325\textwidth}
    \centering
    \includegraphics[width=\linewidth, height=4.6cm]{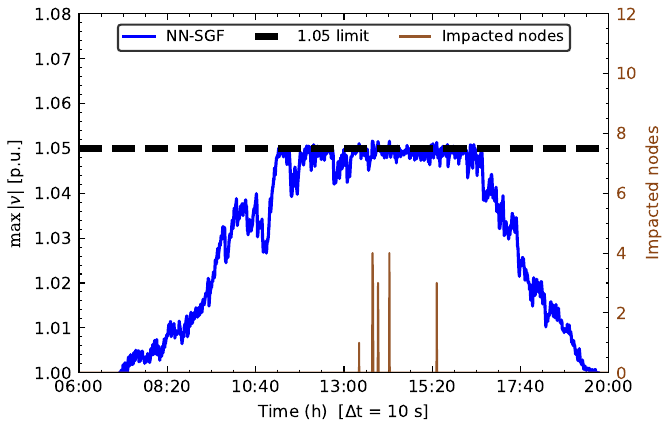}
    \caption{Online NN-SGF in Algorithm~\ref{alg:nn-online}.}
\end{subfigure}
\hfill
\begin{subfigure}[t]{0.325\textwidth}
    \centering
    \includegraphics[width=\linewidth, height=4.6cm]{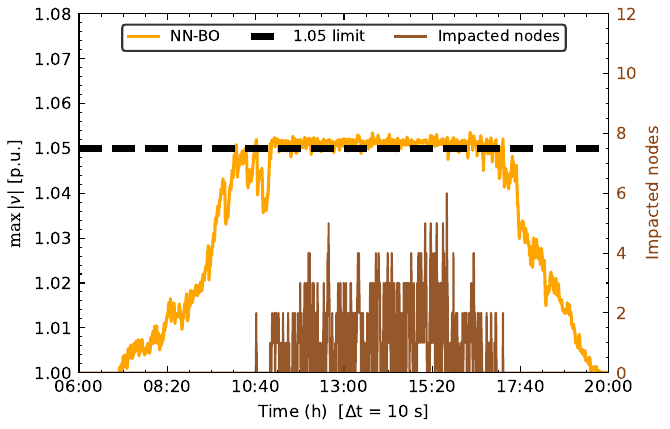}
    \caption{Neural network trained to emulate  $\textsf{U}^{\textsf{lm}}(\vecs_l,\vectheta)$.}
\end{subfigure}
\begin{subfigure}[t]{0.52\textwidth}
    \centering
    \includegraphics[height=3.8cm]{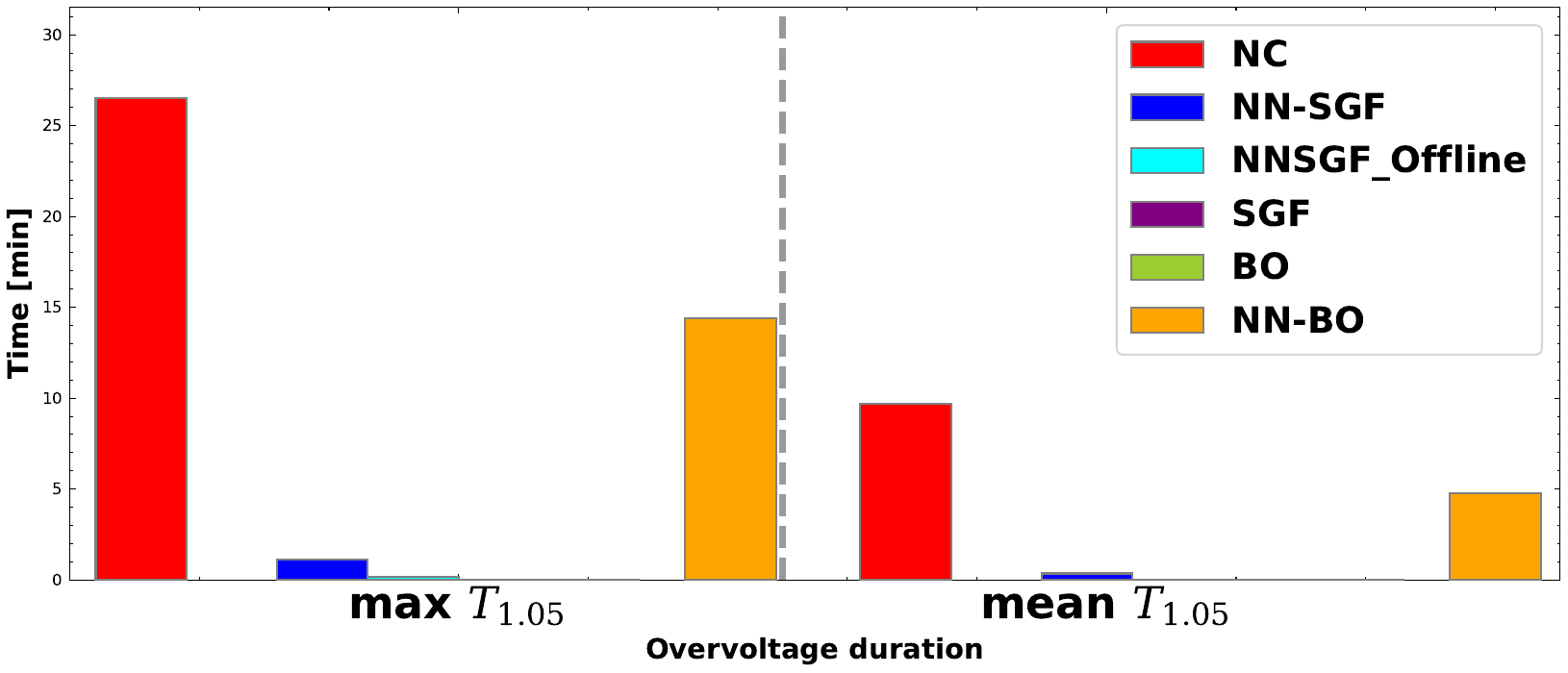}
    \caption{Overvoltage duration.}
\end{subfigure}
\caption{Highest voltage profile  and number of nodes experiencing overvoltages with different optimization methods: \emph{(a)} batch optimization (BO) using IPOPT;  \emph{(b)}  the proposed NN-SGF~\eqref{eq:main_NN}, implemented in a batch configuration as in Figure~\ref{fig:proposed_arch}(center);
\emph{(c)} the proposed NN-SGF~\eqref{eq:main_NN}, implemented in an online feedback configuration as in Figure~\ref{fig:proposed_arch}(left); \emph{(d)} method where a neural network is trained to emulate BO solutions (NN-BO).  Overvoltage durations for the different strategies are compared in \emph{(e)};  we consider the no control (NC) setup as in Figure~\ref{fig:voltage_nc} and the SGF.    }
\label{fig:voltage_comparison}
\vspace{-.5cm}
\end{figure*}


\begin{table}[t!]
\centering
\caption{Training of NN-SGF  (s2) and NN-BO (s3)  }
\label{tab:training}
\begin{tabular}{l|cc}
\hline
\textbf{Method} & \textbf{Training points} & \textbf{Mean Squared Error test}   \\
\hline
NN-SGF (s2) & 60,000 & 1.7 $\times 10^{-6}$  \\
NN-BO (s3) & 80,000 & 8.3 $\times 10^{-5}$  \\
\end{tabular}
\end{table}

Figure~\ref{fig:voltage_comparison}  illustrate how different strategies perform in regulating voltage magnitudes within the bounds $[0.95, 1.05]$~p.u.; Figure~\ref{fig:voltage_comparison} considers the \emph{maximum voltage profile} across the system at every time step, as well as the \emph{number of nodes that experience overvoltages}. {These metrics are evaluated over the monitored set $\mathcal{M}=\{93,54,85,84,48,39,91,71,62,67\}$, selected to cover critical voltage locations, such as end-of-feeder buses and buses that approach or violate voltage limits in the baseline no-control case (Fig.~3); we report $\max_{j\in\mathcal{M}}|v_j(t)|$ (left axis) together with $N_{\mathrm{imp}}(t):=\bigl|\{\,j\in\mathcal{M}:\ |v_j(t)|>\bar V\,\}\bigr|$ (secondary axis).} The proposed NN–SGF method maintains voltages tightly bounded across all monitored nodes, with only brief and mild excursions slightly above 1.05~p.u. These short-duration deviations are well within the tolerance accepted by distribution utilities and do not compromise protection schemes. The BO method, which solves the full AC OPF offline to convergence, is used as a benchmark. The NN-BO approach, trained to emulate the BO setpoints directly, exhibits significantly more overvoltage excursions than the NN-SGF, reflected in its higher $\max T_{1.05}$ and $\mathrm{mean}\ T_{1.05}$ as shown in~\ref{fig:voltage_comparison}(e). This confirms that  approaches that attempts to learn solutions to the OPF directly cannot ensure feasibility. Importantly, NN-SGF delivers effective online voltage regulation without any iterative optimization and even outperforms widely used schemes such as Volt/Var Control (VVC) and online primal dual methods investigated in~\cite{colot2024optimal}, which typically suffer from slower response or larger transient violations. Overall, NN-SGF strikes the best balance among voltage compliance, computational efficiency, and system protection, all without introducing operational concerns.


\subsection{Computational times}

We assess the computational time of the proposed method.  In Table~\ref{tab:runtime_comparison}, we first  consider an online  implementation of the SGF and of the NN-SGF. We recall that the ``Online step'' refers to the setup in Figure~\ref{fig:proposed_arch}(left) where the one evaluation of the SGF (resp., the NN-SGF) is used to generate new setpoints $\vecu(t)$, and the setpoints are sent to the inverters. We averaged the runtime over the full simulation horizon \( t \in [\text{06{:}00}, \text{20{:}00}] \).
Obviously, the computational time of the NN-SGF is much lower, as the SGF involves solving the constrained QP in~\eqref{eq:controller_approximated} to obtain the setpoint \( \dot{\vecu}(t_n) = \eta F_{\textsf{in}}(\vecu(t_n), \vecxi(t_n), \vectheta(t_n)) \). 

As a point of comparison, we consider the average computational time required by \texttt{IPOPT} to solve the AC OPF for the network in Figure~\ref{fig:network}, which is reported is Table~\ref{tab:runtime_comparison}.
 
Overall, an \emph{online} implementation of the NN-SGF achieves a \( \sim\!297\times \) speedup over BO and a \( \sim\!45\times \) speedup over the online version of the SGF proposed in~\cite{colot2024optimal}.  

\begin{table}[h]
\centering
\caption{Average computation times (in seconds). The times do not include the delay in measuring voltages (for SGF) or loads (for BO).  }
\label{tab:runtime_comparison}
\begin{tabular}{l|ccc}
\hline
\textbf{Method} & \textbf{Online step} & \textbf{Offline implem.} &  \textbf{Offline} \\
&  \textbf{Fig.~\ref{fig:proposed_arch}(left)} & \textbf{Fig.~\ref{fig:proposed_arch}(center)} & \textbf{solution} \\
\hline
SGF   &  0.1158 & 1.181 & --  \\
NN-SGF & 0.0026 &  0.236 & -- \\
BO (\texttt{IPOPT}) & -- &  -- & 0.771 \\
NN-BO & -- & -- & 0.0021 \\
\hline
\end{tabular}
\end{table}

We note that, even when NN–BO is trained using consistent high-voltage local minimizers, the resulting input–output map can be difficult to approximate robustly when deployed in a one-shot, open-loop fashion. The observed performance gap between NN–SGF and NN–BO is therefore not primarily due to set-valuedness in the sense of~\cite{sun2018learning}, but rather to how feasibility and approximation errors are handled. In particular, NN–BO applies a learned OPF solution open loop, so that small approximation errors—especially near active constraints—can directly induce voltage violations. By contrast, NN–SGF learns a safe update direction from a strongly convex QP and embeds it in feedback, enforcing (practical) feasibility by construction through forward invariance and explicitly accounting for approximation errors, which leads to improved robustness near constraint boundaries.

We also consider the offline implementation. Here, the ``Offline solution'' refers to Figure~\ref{fig:proposed_arch}(center), where each iteration involves one evaluation of the SGF (resp., the NN-SGF) and one solution to the PF. Again, the PF equations are solved  using \texttt{pandapower}, and the average execution time of \texttt{pandapower} was 0.021 seconds. On average, the proposed scheme implemented in an offline fashion required less than 10 iterations, yielding the upper bounds provided in Table~\ref{tab:runtime_comparison}. 
We note that the proposed online NN-SGF method requires only voltage feedback measurements at buses in the monitored set $\mathcal{M}$. In contrast, BO and NN-BO rely on system-wide measurements or estimates of all non-controllable loads in the network. In practice, acquiring such system-wide information typically requires distribution-level state estimation, which may introduce non-negligible latency beyond the solver runtime. As a result, the end-to-end execution time of BO and NN-BO can be substantially larger in practice than that of NN-SGF~\cite{dehghanpour2018survey,muscas2014effects,angioni2015impact}.

\subsection{{Cumulative cost}}

\begin{figure}[t]
    \centering
    \includegraphics[width=.95\linewidth]{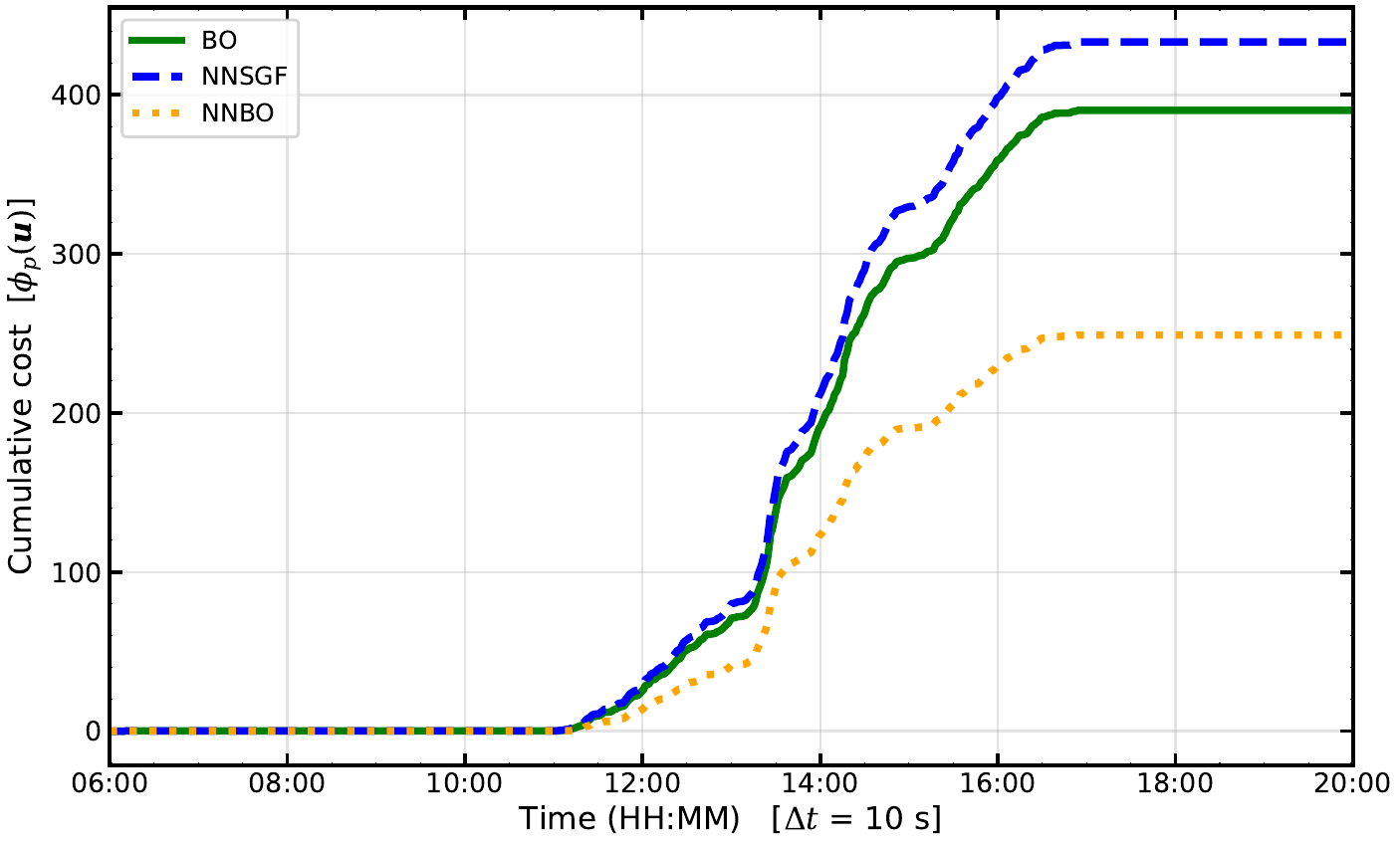}
    \caption{{Cumulative cost over time for BO, NN--SGF, and NN--BO, where the cost penalizes active-power curtailment and reactive-power utilization.}}
    \label{fig:cumcost}
    \vspace{-.4cm}
\end{figure}

Figure~\ref{fig:cumcost} reports the cumulative cost associated with active-power curtailment and reactive-power utilization. Among the considered methods, NN-BO attains the lowest cumulative cost, primarily because it performs less active-power curtailment—thus exporting more available DER generation—and applies smaller reactive-power adjustments, corresponding to weaker corrective actions. However, this apparent advantage comes at the expense of voltage regulation performance. In particular, NN-BO frequently operates outside the admissible voltage region and therefore avoids the curtailment. 

In contrast, BO and NN-SGF incur higher cumulative costs because they actively enforce operational constraints. When voltages approach or violate prescribed limits, these methods increase corrective actions through active-power curtailment and/or reactive-power support to drive the system back within allowable bounds. Importantly, NN-SGF enforces feasibility online via feedback-based updates and achieves rapid voltage recovery under fast-changing operating conditions. 

\section{Theoretical Analysis}
\label{sec:convergence}

In this section, we analyze the convergence and the ability to generate feasible points of our proposed method~\eqref{eq:main_NN}. We start with the following assumption, which imposes some mild regularity assumptions on a neighborhood of a strict locally optimal solution of the AC OPF. 

\begin{assumption}[\textit{Regularity conditions}]
\label{as:openloop}
\normalfont    Assume that the problem~\eqref{eq:AC-OPF} is feasible and  let $\vecu^*$ be an isolated KKT
    point for~\eqref{eq:AC-OPF} and a local minimizer, for given $\vecpl,\vecql$. Assume that: 

    \noindent \textit{i)}
    Strict complementarity slackness~\cite{fiacco1976sensitivity} and the linear independence constraint qualification (LICQ)~\cite{hauswirth2018generic} hold at $ \vecu^*$.
    
    \noindent \textit{ii)}
    The maps $\vecu \mapsto \phi_p(\vecu)$, $\vecu \mapsto \phi_v(|\vecv(\vecu;\vecpl,\vecql)|) $,  $\vecu\mapsto |\vecv(\vecu;\vecpl,\vecql)|$, and $\vecu\mapsto |\veci(\vecu;\vecpl,\vecql)|$ are twice continuously differentiable over $\mathcal{B}(\vecu^*,r_1)$, and their Hessian matrices are positive semi-definite at $\vecu^*$.
    
    \noindent \textit{iii)}
    The Hessian $\nabla^2  \phi_p(\vecu^*)$ is positive definite.
     \hfill $\Box$
\end{assumption}

This assumption is supported by the results of~\cite{hauswirth2018generic} and used in~\cite{colot2024optimal}. We also impose the following assumptions on the approximation and training errors. 

\begin{assumption}[\textit{Approximation errors}] 
\label{as:bound_approx_error}
\normalfont
    $\exists ~ E_v < + \infty, E_{J_v} < + \infty$ such that $\||\vecv(\vecu; \vecs_l)| - (\boldsymbol{\Gamma}_v \vecu  +  \bar{\vecv}(\vecs_l)) \|_\infty \leq E_v$,  and $\|\boldsymbol{\Gamma}_v -J_v(\vecu,\vecs_l)  \| \leq E_{J_v}$ for any  $\vecu \in \mathcal{B}(\vecu^*,r_1)$. \hfill $\Box$
\end{assumption}

\begin{assumption}[\textit{Measurement errors}] 
\label{as:bound_meas_error}
\normalfont
$\exists ~ \epsilon_n < + \infty$ such that $\|\vecn\|_\infty \leq \epsilon_n$. Consequently, we let $\epsilon_v < + \infty$ and $~ \epsilon_i < + \infty$ such that $\|\vecn_v\|_\infty \leq \epsilon_v$ and $\|\vecn_i\|_\infty \leq \epsilon_i$, respectively. \hfill $\Box$
\end{assumption}

\begin{assumption}[\textit{Training errors}] 
\label{as:bound_train_error}
\normalfont
$\exists ~ E^{\textsf{NN}} < + \infty$ such that $\|\mathcal{F}^{\textsf{NN}}(\vecu, \vecxi, \vectheta) - F_{\text{ln}}(\vecu, \vecxi, \vectheta\| \leq E^{\textsf{NN}}$  for all $(\vecu, \vecxi, \vectheta) \in \mathcal{C}_{\text{train}} \times \mathcal{E}_{\text{train}} \times \Theta_{\text{train}}$. 
 \hfill $\Box$
\end{assumption}

Since the line currents $\veci$ can be computed from $\vecv$ via Ohm's Law, Assumption~\ref{as:bound_approx_error}
implies that $\exists ~ E_i < + \infty, E_{J_i} < + \infty$ such that $\||\veci(\vecu; \vecs_l)| - (\boldsymbol{\Gamma}_i \vecu  +  \bar{\veci}(\vecs_l)) \|_\infty \leq E_i$  and $\|\boldsymbol{\Gamma}_i -J_i(\vecu,\vecs_l)  \| \leq E_{J_i}$ for any  $\vecu \in \mathcal{B}(\vecu^*,r_1)$. Assumptions~\ref{as:bound_approx_error}-\ref{as:bound_meas_error} are motivated by the fact that the error of the linear approximation is small~\cite{bolognani2015existence,sadnan2021learning}, and that in realistic monitoring and SCADA systems, the measurement of the voltage magnitudes are affected by a negligible error~\cite{angioni2015impact}. 

Lastly, Assumption~\ref{as:bound_train_error} follows from the approximation capabilities of neural networks over compact sets~\cite{hornik1991approximation,duan2023minimum}. 

To proceed, denote as $\boldsymbol{\Psi}_v :=  \boldsymbol{\Gamma}_v - \boldsymbol{J}_v(\vecu;\vecs_l)$ and $\boldsymbol{\Psi}_i :=  \boldsymbol{\Gamma}_i -  \boldsymbol{J}_i(\vecu;\vecs_l)$ the errors in the computation of the Jacobian, and define the sets $\mathcal{E}_v = \{\boldsymbol{\Psi}_v: \|\boldsymbol{\Psi}_v\| \leq E_{J_v}\}$ and $\mathcal{E}_i = \{\boldsymbol{\Psi}_i: \|\boldsymbol{\Psi}_i\| \leq E_{J_i}\}$. Let $\mathcal{E}_n := \{\vecn: \|\vecn\| \leq \epsilon_n\}$. 
Define the map $F_{\textsf{m}}(\vecu, \vecn, \boldsymbol{\Psi}_v, \boldsymbol{\Psi}_i)$ as
\begin{align}
    &  F_{\textsf{m}}(\vecu, \vecn, \boldsymbol{\Psi}_v, \boldsymbol{\Psi}_i) \label{eq:Fm} \\
    & := \arg \min_{\vecz \in \mathbb{R}^{2G}}  \|\vecz  + \nabla \phi_p(\vecu) + (\boldsymbol{J}_v(\vecu;\vecs_l) + \boldsymbol{\Psi}_v)^\top\nabla \phi_v(\vecV)\|^2 \nonumber  \\
    & \hspace{1cm}  \textrm{s.t.}  -(\boldsymbol{J}_v(\vecu;\vecs_l) + \boldsymbol{\Psi}_v)^\top \vecz  \leq -\beta \left(\mathbf{1} \underline{V}- (|\vecv| + \vecn_v) \right) \nonumber \\
    & \hspace{1.8cm} (\boldsymbol{J}_v(\vecu;\vecs_l) + \boldsymbol{\Psi}_v)^\top \vecz  \leq -\beta \left((|\vecv| + \vecn_v) -\bar{V} \mathbf{1}\right) \nonumber \\  & \hspace{1.8cm} (\boldsymbol{J}_i(\vecu;\vecs_l) + \boldsymbol{\Psi}_i)^\top \vecz \leq -\beta \left((|\veci| + \vecn_i) -\bar{I} \mathbf{1}\right) \nonumber \\
    &  \hspace{1.7cm} \boldsymbol{J}_{\ell_i}(\vecu_i)^\top \vecz  \leq - \beta \ell_i(\vecu), ~ i \in \mathcal{G} \nonumber 
\end{align}
which is a representation of $F_{\textsf{lm}}(\vecu, \vecxi,\vectheta)$
emphasizing the dependence on the errors;  note also that $F(\vecu, \vecxi,\vectheta) = F_{\textsf{m}}(\vecu, \boldsymbol{0}, \boldsymbol{0}, \boldsymbol{0})$. With this notation, we assume the following.

\begin{assumption}[\textit{MFCQ and constant-rank conditions}] 
\label{as:Regularity}
\normalfont For any $\vecu  \in \mathcal{B}(\vecu^*,r_1)$, and any $\boldsymbol{\Psi}_v$, $\boldsymbol{\Psi}_i$, and $\boldsymbol{n}_v$, $\boldsymbol{n}_i$ satisfying Assumptions~\ref{as:bound_approx_error}-\ref{as:bound_meas_error}, the problem~\eqref{eq:controller_approximated} is feasible,  and satisfies the Mangasarian-Fromovitz Constraint Qualification  and the constant-rank condition~\cite{liu1995sensitivity}. 
 \hfill $\Box$
\end{assumption}

Next, we present the following intermediate result. 

\begin{lemma}[\textit{Continuity of $F_{\textsf{m}}(\vecu, \vecn, \boldsymbol{\Psi}_v, \boldsymbol{\Psi}_i)$}] 
\label{lem:lipschitz}
\normalfont

Let Assumption~\ref{as:Regularity} hold, and assume that $\vecu \mapsto \phi_p(\vecu)$, $\vecV \mapsto \phi_v(\vecV) $ are twice continuously differentiable over $\mathcal{B}(\vecu^*,r_1)$ and for any $\vecV$ Then: 

\noindent (i) For any $\vecn \in \mathcal{E}_n$, $\boldsymbol{\Psi}_v \in \mathcal{E}_v$, and $\boldsymbol{\Psi}_i \in \mathcal{E}_i$, $\vecu \mapsto F_{\textsf{m}}(\vecu, \vecn, \boldsymbol{\Psi}_v, \boldsymbol{\Psi}_i)$ is locally Lipschitz at $\vecu$, $\vecu \in \mathcal{B}(\vecu^*,r_1)$. 

\noindent (ii) For any $\vecu \in \mathcal{B}(\vecu^*,r_1)$, $\boldsymbol{\Psi}_v \in \mathcal{E}_v$, and $\boldsymbol{\Psi}_i \in \mathcal{E}_i$, $\vecn \mapsto F_{\textsf{m}}(\vecu, \vecn, \boldsymbol{\Psi}_v, \boldsymbol{\Psi}_i)$ is  Lipschitz with constant $\ell_n \geq 0$  over $\mathcal{E}_n$.

\noindent (iii) For any $\vecu \in \mathcal{B}(\vecu^*,r_1)$, $\vecn \in \mathcal{E}_n$, $\boldsymbol{\Psi}_i \in \mathcal{E}_i$, $\boldsymbol{\Psi}_v \mapsto F_{\textsf{m}}(\vecu, \vecn, \boldsymbol{\Psi}_v, \boldsymbol{\Psi}_i)$ is  $\ell_{J_{v}}$-Lipschitz over $\mathcal{E}_v$. 

\noindent (iv) For any $\vecu \in \mathcal{B}(\vecu^*,r_1)$, $\vecn \in \mathcal{E}_n$, $\boldsymbol{\Psi}_v \in \mathcal{E}_v$, $\boldsymbol{\Psi}_i \mapsto F_{\textsf{m}}(\vecu, \vecn, \boldsymbol{\Psi}_v, \boldsymbol{\Psi}_i)$ is  $\ell_{J_{i}}$-Lipschitz  over $\mathcal{E}_i$. \hfill $\Box$
\end{lemma}

Lemma~\ref{lem:lipschitz} follows from~\cite[Theorem 3.6]{liu1995sensitivity}, and by the compactness of the sets $\mathcal{E}_n$,  $\mathcal{E}_{v}$ and $\mathcal{E}_i$. 
To state the main result, recall that $\vecu^*$ is a local optimizer of \eqref{eq:AC-OPF}. Define $\vecxi^*:=H(\vecu^*;\vecs_l)$, 
$\boldsymbol{J}_F :=\frac{\partial F (\vecu,H(\vecu; \vecs_l),\vectheta)}{\partial \vecu}\mid_{\vecu=\vecu^*}$, 
$e_1 :=-\lambda_{\max}(\boldsymbol{J}_F)$, 
and 
$e_2 :=-\lambda_{\min}(\boldsymbol{J}_F)$. 
Then, from~\cite{allibhoy2023control}, we can write the dynamics as $F (\vecu,H(\vecu; \vecs_l),\vectheta) = \boldsymbol{J}_F (\vecu-\vecu^*) + g(\vecu)$, where $g$ satisfies $\|g(\vecu)\|\leq L \|\vecu-\vecu^*\|^2$, $\forall \vecu\in \mathcal{B}(\vecu^*,r_2)$, for some $L> 0$ and $r_2>0$. Define $r:=\min\{r_1,r_2\}$ and $s_{\min}$ as: $s_{\min} = 0$ if $ r\geq \frac{e_1}{L}$, and $s_{\min} = 1-\frac{r L}{e_1}$ if $r\geq \frac{e_1}{L}$. We are ready to state the main result.  

\begin{theorem}[\emph{Stability}] 
\label{thm:stability}
\normalfont
 Consider the OPF problem~\eqref{eq:AC-OPF} satisfying  Assumption \ref{as:steadyStateMap}. Let Assumptions~\ref{as:bound_approx_error}--\ref{as:bound_train_error} hold for the linear model and the training, and let Assumption~\ref{as:Regularity} hold for~\eqref{eq:controller_approximated}. Let $\vecu(t)$, $t \geq t_0$,  be the unique trajectory of~\eqref{eq:main_NN}. Let $\epsilon := \ell_{J_v} E_{J_v} + \ell_{J_i} E_{J_i} + \ell_{n} \epsilon_n +  \epsilon^{\textsf{NN}}$, and assume that the set $ \mathcal{R}:=\left\{s:s_{\min}<s\leq 1,~e_1^{-3} e_2 L \epsilon <s-s^2\right\}$ is not empty. Then, for any $s\in\mathcal{R}$, it holds that 
\begin{align}
 &\|\vecu(t)-\vecu^*\|
 \leq \sqrt{\frac{e_2}{e_1}} e^{- e_1\eta s(t-t_0)}\|\vecu(t_0)-\vecu^*\|  + \frac{e_2 \epsilon}{ s e_1^2}  \label{eq:transient_u},
\end{align}
for any $\vecu(t_0)$ such that $\|\vecu(t_0)-\vecu^*\|\leq \sqrt{\frac{e_1}{e_2}}\frac{e_1}{L}(1-s)$. 
\hfill $\triangle$
\end{theorem}

From Theorem~\ref{thm:stability}, one can see that the asymptotic error can be reduced by improving the approximation accuracy of the neural network (i.e., reducing $\epsilon^{\textsf{NN}}$) or by improving the linear approximation (i.e., reducing $E_{J_v}$ and $E_{J_i}$). Practically,  the errors in the voltages and currents  (i.e.,  $\epsilon_n$) are negligible.  As a technical detail, the requirement that  $\mathcal{R}$ is not empty guarantees that $\vecu(t)$  never exits the region of attraction of $\vecu^*$.

The following result characterizes the feasibility of $\vecu(t)$.  

\begin{proposition}[\textit{Practical feasibility}] 
\label{prop:invariance}
\normalfont
Let  the conditions in Theorem~\ref{thm:stability} be satisfied, and let $\vecu(t)$, $t \geq t_0$,  be the unique trajectory of~\eqref{eq:main_NN}. Define the set 
$\mathcal{S}_{e} := \mathcal{S}_{e,v} \cap \mathcal{S}_{e,i}$,
\begin{subequations}
\label{eq:setS}   
\begin{align}
    \mathcal{S}_{e,v} & := \{\vecu \in \mathcal{C}: \underline{V}_e \leq |v_j(\vecu; \vecs_l)| \leq \overline{V}_e , j \in \mathcal{M}\} \\
    \mathcal{S}_{e,i} & := \{\vecu \in \mathcal{C}: |i_i(\vecu; \vecs_l)| \leq \overline{I}_e, j \in \mathcal{L}\} 
\end{align}
\end{subequations}
where $\underline{V}_e = \underline{V} - \epsilon_v - 2 E_v - \beta^{-1}\|\boldsymbol{\Gamma}_v\| \epsilon^{\textsf{NN}}$, $\overline{V}_e = \overline{V} + \epsilon_v + 2 E_v + \beta^{-1} \|\boldsymbol{\Gamma}_v\| \epsilon^{\textsf{NN}}$, and $\overline{I}_e = \overline{I} + \epsilon_i + 2 E_i + \beta^{-1} \|\boldsymbol{\Gamma}_i\| \epsilon^{\textsf{NN}}$. Then, the neural network-based algorithm~\eqref{eq:main_NN} renders a set $ \mathcal{S}_s$, with 
$\mathcal{S} \subseteq  \mathcal{S}_s \subseteq \mathcal{S}_e$, forward invariant.
\hfill $\triangle$
\end{proposition}

We recall that $E_v$ and $E_i$ are defined in Assumption~\ref{as:bound_approx_error}. The set $\mathcal{S}_s$ is an auxiliary set defined through linearized voltage and current magnitudes and is instrumental for establishing forward invariance of the NN-SGF dynamics. The formal definition of $\mathcal{S}_s$ is provided in the appendix. 
The set $\mathcal{S}_e$, instead, is expressed in terms of the actual voltage and current magnitudes and contains $\mathcal{S}_s$ by accounting for the linearization and neural-network approximation errors. This construction leads to the factor \(2\) appearing in the definitions of 
\(\underline V_e\), \(\overline V_e\), and \(\overline I_e\).
We note that $\mathcal{S}_s$ in  Proposition~\ref{prop:invariance} is an inflation of the set of  feasible voltages and currents $\mathcal{S}$ specified in the AC OPF. Hence the terminology ``practical feasibility''.  Indeed, when these errors are small, the constraint violation is practically negligible. The proofs of all results are provided in the appendix.

We provide the following remarks. 

\begin{remark}[\emph{Practical forward invariance}]
{\rm If $\vecu(t_0)\in\mathcal{S}_s$, then $\vecu(t)\in\mathcal{S}_s$ for all $t\ge t_0$. Consequently, the control input generated by~\eqref{eq:main_NN} is practically feasible both for the online implementation shown in Fig.~\ref{fig:proposed_arch}(left) and when the offline procedure in Fig.~\ref{fig:proposed_arch}(center) is terminated prior to convergence. Moreover, the bounds $\underline{V}_e$, $\overline{V}_e$, and $\overline{I}_e$ provide an explicit quantification of the sets of voltage magnitudes and currents that are achievable under~\eqref{eq:main_NN} as a function of the bounds $\epsilon_v$, $\epsilon_i$, $E_v$, $E_i$, and the neural-network approximation error $\epsilon^{\textsf{NN}}$. For instance, in the experiments reported in Section~\ref{sec:results}, the maximum value of $\|\boldsymbol{\Gamma}_v\|$ and $\|\boldsymbol{\Gamma}_i\|$ was $0.0982$, $E_v \approx 10^{-2}$ p.u. and $\epsilon^{\textsf{NN}} \approx 1.7\times 10^{-6}$. The induced constraint violations were practically negligible.} \hfill $\Box$
\end{remark}

\begin{figure}[t]
    \centering
        \includegraphics[width=6.0cm]{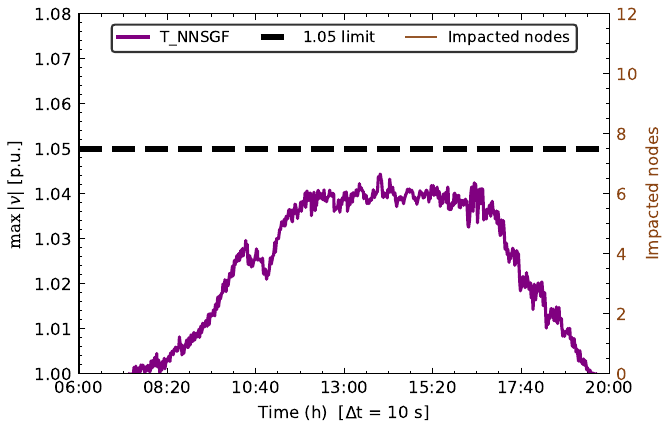}
    \caption{Numerical experiment with voltage constraints tightened based on the result of Proposition~\ref{prop:invariance}. The NN-SGF is trained based on  lower and upper voltage constraints given by $\underline{V} + \epsilon_v + 2 E_v + \beta^{-1}\|\boldsymbol{\Gamma}_v\| \epsilon^{\textsf{NN}}$ and $\overline{V} - \epsilon_v - 2 E_v - \beta^{-1} \|\boldsymbol{\Gamma}_v\| \epsilon^{\textsf{NN}}$, respectively. One can see that the voltage limits are satisfied strictly.}
    \label{fig:tightened_sgf}
\end{figure}

\begin{remark}[\emph{OPF with constraint tightening}]
{\rm When bounds on $\epsilon_v$, $\epsilon_i$, $E_v$, $E_i$, and the neural-network approximation error $\epsilon^{\textsf{NN}}$ are available or can be estimated numerically, the constraints of the original AC OPF problem~\eqref{eq:AC-OPF} can be appropriately tightened so that the control policy generated by~\eqref{eq:main_NN} renders the feasible set $\mathcal{S}$ forward invariant. In particular, leveraging the result of Proposition~\ref{prop:invariance}, such a tightening guarantees that voltage and current limits are satisfied despite modeling, linearization, and approximation errors. This can be achieved by training the proposed NN-SGF using tightened voltage bounds given by $\underline{V} + \epsilon_v + 2 E_v + \beta^{-1}\|\boldsymbol{\Gamma}_v\| \epsilon^{\textsf{NN}}$ and $\overline{V} - \epsilon_v - 2 E_v - \beta^{-1} \|\boldsymbol{\Gamma}_v\| \epsilon^{\textsf{NN}}$, respectively.

We evaluated the performance of the NN-SGF under such tightened constraints. In the numerical experiments, we assumed negligible measurement noise (i.e., $\epsilon_v=0$), estimated the voltage linearization error bound as $E_v \approx 10^{-2}$ p.u., and obtained a neural-network approximation error of $\epsilon^{\textsf{NN}} = 1.7\times 10^{-6}$. Substituting these values into Proposition~\ref{prop:invariance} yields tightened voltage limits of $0.9604$ and $1.0396$. The resulting AC OPF problem with tightened constraints remained feasible for all tested load and PV generation scenarios. The voltage profiles obtained in this experiment are reported in Fig.~\ref{fig:tightened_sgf}. As observed, the voltage magnitude constraints are never violated, demonstrating a robust mechanism to account for training and modeling errors. As expected, the resulting solution is conservative, owing to the bounding techniques employed in the derivation of Proposition~\ref{prop:invariance} (e.g., repeated use of triangle inequalities, as shown in the appendix). \hfill $\Box$
}
\end{remark}

\section{Conclusions}
\label{sec:conclusions}
We have proposed a solution method for solving the AC OPF where a neural network is used to approximate the solution of a convex QP defining the safe gradient flow. Our approach is shown to lead to both feedback-based online  implementations and offline solutions based on power flow computations. Compared to existing methods that rely on neural networks, our algorithm ensures that the DERs' setpoints are practically feasible and it ensures convergence to a neighborhood of a strict local optimizer of the AC OPF. These guarantees are important for power systems optimization tasks, as operating limits must be satisfied for a safe power delivery. Results are derived based on  assumptions that are sufficient and conservative by design, as they rely on worst-case bounds to guarantee feasibility and forward invariance in the presence of modeling, measurement, and learning errors. In practice, the corresponding error bounds can be estimated empirically and were  small in the numerical results. Future work will consider active set learning approaches for solving the QP associated with the SGF in the context of AC OPF in distribution grids~\cite{deka2019learning}, and will formalize results in terms of asymptotic stability under a discrete-time version of 
$\mathcal{F}^{\textsf{NN}}\!\big(\boldsymbol{u},\boldsymbol{\xi},\boldsymbol{\theta}\big)$, arising from, for example, forward-Euler methods. Finally, we will explore the derivation of formal results in terms of tracking and feasibility when the AC OPF problem is time-varying, under appropriate  assumptions on the variability of the problem~\cite{tang2022running}.

\bibliographystyle{ieeetr} 
\bibliography{references}

\appendices
\section*{APPENDIX}
\label{app:robustness}

\subsection{Proof of Theorem~\ref{thm:stability}}

Recall that \( \vecV := |\vecv(\vecu; \vecs_l)| + \vecn_v\), \( \vecI := |\veci(\vecu; \vecs_l)| + \vecn_i\), $\vecxi = (\vecV, \vecI)$; to streamline notation, we will use $|\vecv|$ and $|\veci|$ to denote the error-free measurements or computations of voltage magnitudes and currents magnitudes. Recall that \( F(\vecu, \vecxi, \vectheta) = F_{\textsf{m}}(\vecu, \boldsymbol{0}, \boldsymbol{0}, \boldsymbol{0}) \). We express the NN-SGF controller as \( \dot{\vecu} = \mathcal{F}^{\textsf{NN}}(\vecu, \vecxi, \vectheta) \), where \( \vectheta = (\vectheta_{u,1}, \ldots, \vectheta_{u,G}, \underline{V}, \overline{V}, \overline{I}) \) contains the constraint parameters of the AC OPF. Rewrite the NN-SGF controller as:
\begin{align*}
\dot{\vecu} &= \eta\, \mathcal{F}^{\textsf{NN}}(\vecu, \vecxi, \vectheta) \\
&= \eta\, \underbrace{F_{\textsf{m}}(\vecu, \boldsymbol{0}, \boldsymbol{0}, \boldsymbol{0})}_{\text{nominal}} \\
&\quad + \eta\, \underbrace{\left[ F_{\textsf{m}}(\vecu, \vecn, \boldsymbol{\Psi}_v, \boldsymbol{\Psi}_i) - F_{\textsf{m}}(\vecu, \vecn, \boldsymbol{0}, \boldsymbol{0}) \right]}_{\text{Jacobian error}} \\
&\quad + \eta\, \underbrace{\left[ F_{\textsf{m}}(\vecu, \vecn, \boldsymbol{0}, \boldsymbol{0}) - F_{\textsf{m}}(\vecu, \boldsymbol{0}, \boldsymbol{0}, \boldsymbol{0}) \right]}_{\text{measurement error}} \\
&\quad + \eta\, \underbrace{\left[ \mathcal{F}^{\textsf{NN}}(\vecu, \vecxi, \vectheta) - F_{\textsf{ln}}(\vecu, \vecxi, \vectheta) \right]}_{\text{NN training error}}
\end{align*}
where we stress that $F_{\textsf{m}}(\vecu, \boldsymbol{0}, \boldsymbol{0}, \boldsymbol{0}) = F(\vecu, (|\vecv|,|\veci|), \vectheta)$ is the nominal controller~\eqref{eq:controller_ideal}.The NN-SGF controller is thus interpreted as a perturbation of the nominal gradient flow \( F_{\textsf{m}}(\vecu, \boldsymbol{0}, \boldsymbol{0}, \boldsymbol{0}) \), which is differentiable at the strict local minimizer \( \vecu^* \) (Assumption~\ref{as:openloop}); its Jacobian is defined as
\[
\boldsymbol{J}_F := \left. \frac{\partial F_{\textsf{m}}(\vecu, \boldsymbol{0}, \boldsymbol{0}, \boldsymbol{0})}{\partial \vecu} \right|_{\vecu = \vecu^*},
\]
and is negative definite. {This follows from properties of the safe gradient flow at strict local minimizers; see \cite[Lemma~5.12]{allibhoy2023control} and \cite[Proof of Theorem~5.7(iii)]{allibhoy2023control}. {In particular, all eigenvalues of $\boldsymbol{J}_F$ are strictly negative, so $e_1 := -\lambda_{\max}(\boldsymbol{J}_F)$ and $e_2 := -\lambda_{\min}(\boldsymbol{J}_F)$ are well-defined positive real scalars.}}

Let \( e_1 := -\lambda_{\max}(\boldsymbol{J}_F) \), \( e_2 := -\lambda_{\min}(\boldsymbol{J}_F) \), and define the matrix
\[
P := \int_0^\infty e^{\boldsymbol{J}_F^\top \zeta} e^{\boldsymbol{J}_F \zeta} \, d\zeta,
\]
which satisfies the Lyapunov equation \( P \boldsymbol{J}_F + \boldsymbol{J}_F^\top P = -I \). The matrix \( P \) satisfies the bounds:
\[
\frac{1}{2e_2} \|\vecu - \vecu^*\|^2 \le (\vecu - \vecu^*)^\top P (\vecu - \vecu^*) \le \frac{1}{2e_1} \|\vecu - \vecu^*\|^2.
\]

Define the Lyapunov function
\[
V_1(\vecu) := (\vecu - \vecu^*)^\top P (\vecu - \vecu^*).
\]

We compute:
\begin{align*}
\dot{V}_1(\vecu) 
&= 2(\vecu - \vecu^*)^\top P \dot{\vecu} \\
&= 2\eta (\vecu - \vecu^*)^\top P\, \mathcal{F}^{\textsf{NN}}(\vecu, \vecxi, \vectheta) \\
&= 2\eta (\vecu - \vecu^*)^\top P\, F_{\textsf{m}}(\vecu, \boldsymbol{0}, \boldsymbol{0}, \boldsymbol{0}) \\
&\quad + 2\eta (\vecu - \vecu^*)^\top P\, \left[ F_{\textsf{m}}(\vecu, \vecn, \boldsymbol{\Psi}_v, \boldsymbol{\Psi}_i) - F_{\textsf{m}}(\vecu, \vecn, \boldsymbol{0}, \boldsymbol{0}) \right] \\
&\quad + 2\eta (\vecu - \vecu^*)^\top P\, \left[ F_{\textsf{m}}(\vecu, \vecn, \boldsymbol{0}, \boldsymbol{0}) - F_{\textsf{m}}(\vecu, \boldsymbol{0}, \boldsymbol{0}, \boldsymbol{0}) \right] \\
&\quad + 2\eta (\vecu - \vecu^*)^\top P\, \left[ \mathcal{F}^{\textsf{NN}}(\vecu, \vecxi, \vectheta) - F_{\textsf{ln}}(\vecu, \vecxi, \vectheta) \right].
\end{align*}
Next, we analyze each term. As for the nominal controller, by first-order Taylor expansion~\cite{allibhoy2023control}:
\begin{align*}
F_{\textsf{m}}(\vecu, \boldsymbol{0}, \boldsymbol{0}, \boldsymbol{0}) 
&= F_{\textsf{m}}(\vecu^*, \boldsymbol{0}, \boldsymbol{0}, \boldsymbol{0}) \\
&\quad + \left.\frac{\partial F_{\textsf{m}}(\vecu, \boldsymbol{0}, \boldsymbol{0}, \boldsymbol{0})}{\partial \vecu} \right|_{\vecu = \vecu^*} (\vecu - \vecu^*) 
+ g(\vecu). 
\end{align*}
Additionally, one has that $\| g(\vecu) \| 
\le L \| \vecu - \vecu^* \|^2$, for some $L \geq 0$. Then, $F_{\textsf{m}}(\vecu, \boldsymbol{0}, \boldsymbol{0}, \boldsymbol{0}) = \boldsymbol{J}_F (\vecu - \vecu^*) + \hat{g}(\vecu)$.

The quadratic form evaluates as:
\begin{align*}
&2\eta(\vecu - \vecu^*)^\top P 
F_{\textsf{m}}(\vecu, \boldsymbol{0}, \boldsymbol{0}, \boldsymbol{0}) \\
&= (\vecu - \vecu^*)^\top 
\left(P \boldsymbol{J}_F + \boldsymbol{J}_F^\top P\right) 
(\vecu - \vecu^*) 
+ 2\eta(\vecu - \vecu^*)^\top P \hat{g}(\vecu) .
\end{align*}
Using the Lyapunov identity \( P \boldsymbol{J}_F + \boldsymbol{J}_F^\top P = -I \), the bound \( \|P\| \le \frac{1}{2e_1} \), and 
\( \|g(\vecu)\| \le L \|\vecu - \vecu^*\|^2 \), we conclude:
\[
2\eta (\vecu - \vecu^*)^\top P F_{\textsf{m}}(\vecu, \boldsymbol{0}, \boldsymbol{0}, \boldsymbol{0}) 
\le -\eta \| \vecu - \vecu^* \|^2 + \frac{\eta L}{e_1} \| \vecu - \vecu^* \|^3.
\]

We now focus on the term related to the error in the Jacobian. Using the triangle inequality we get:
\begin{align*}
&\|F_{\textsf{m}}(\vecu, \vecn, \boldsymbol{\Psi}_v, \boldsymbol{\Psi}_i) 
- F_{\textsf{m}}(\vecu, \vecn, \boldsymbol{0}, \boldsymbol{0})\| \\
&\quad = \|F_{\textsf{m}}(\vecu, \vecn, \boldsymbol{\Psi}_v, \boldsymbol{\Psi}_i) 
- F_{\textsf{m}}(\vecu, \vecn, \boldsymbol{0}, \boldsymbol{\Psi}_i) \\
&\qquad + F_{\textsf{m}}(\vecu, \vecn, \boldsymbol{0}, \boldsymbol{\Psi}_i) 
- F_{\textsf{m}}(\vecu, \vecn, \boldsymbol{0}, \boldsymbol{0})\| \\
&\quad \le \|F_{\textsf{m}}(\vecu, \vecn, \boldsymbol{\Psi}_v, \boldsymbol{\Psi}_i) 
- F_{\textsf{m}}(\vecu, \vecn, \boldsymbol{0}, \boldsymbol{\Psi}_i)\| \\
&\qquad + \|F_{\textsf{m}}(\vecu, \vecn, \boldsymbol{0}, \boldsymbol{\Psi}_i) 
- F_{\textsf{m}}(\vecu, \vecn, \boldsymbol{0}, \boldsymbol{0})\|.
\end{align*}

By Lemma~\ref{lem:lipschitz} and Assumption~\ref{as:bound_approx_error}, there exist constants \( \ell_{J_v} \), \( \ell_{J_i} \) such that:
\begin{align*}
\left\|F_{\textsf{m}}(\vecu, \vecn, \boldsymbol{\Psi}_v, \boldsymbol{\Psi}_i) 
- F_{\textsf{m}}(\vecu, \vecn, \boldsymbol{0}, \boldsymbol{\Psi}_i)\right\| 
&\le \ell_{J_v} E_{J_v}, \\
\left\|F_{\textsf{m}}(\vecu, \vecn, \boldsymbol{0}, \boldsymbol{\Psi}_i) 
- F_{\textsf{m}}(\vecu, \vecn, \boldsymbol{0}, \boldsymbol{0})\right\| 
&\le \ell_{J_i} E_{J_i}.
\end{align*}
Hence,
\[
\left\|F_{\textsf{m}}(\vecu, \vecn, \boldsymbol{\Psi}_v, \boldsymbol{\Psi}_i) 
- F_{\textsf{m}}(\vecu, \vecn, \boldsymbol{0}, \boldsymbol{0})\right\| 
\le \ell_{J_v} E_{J_v} + \ell_{J_i} E_{J_i}.
\]
Using the fact that \( \|P\| \le \frac{1}{2e_1} \), one has that:
\begin{align*}
& 2\eta (\vecu - \vecu^*)^\top P 
\left( F_{\textsf{m}}(\vecu, \vecn, \boldsymbol{\Psi}_v, \boldsymbol{\Psi}_i) 
- F_{\textsf{m}}(\vecu, \vecn, \boldsymbol{0}, \boldsymbol{0}) \right) \\
&\quad \le \frac{2\eta}{2e_1} \left\|F_{\textsf{m}}(\vecu, \vecn, \boldsymbol{\Psi}_v, \boldsymbol{\Psi}_i) 
- F_{\textsf{m}}(\vecu, \vecn, \boldsymbol{0}, \boldsymbol{0})\right\| \cdot \|\vecu - \vecu^*\| \\
&\quad = \frac{\eta}{e_1} \left( \ell_{J_v} E_{J_v} + \ell_{J_i} E_{J_i} \right) \cdot \|\vecu - \vecu^*\|.
\end{align*}

Next, by Lemma~\ref{lem:lipschitz} and Assumption~\ref{as:bound_meas_error}, there exists \( \ell_n \) such that
\[
\left\| F_{\textsf{m}}(\vecu, \vecn, \boldsymbol{0}, \boldsymbol{0}) 
- F_{\textsf{m}}(\vecu, \boldsymbol{0}, \boldsymbol{0}, \boldsymbol{0}) \right\| 
\le \ell_n \epsilon_n.
\]
Then:
\begin{align*}
&2\eta (\vecu - \vecu^*)^\top P 
\left( F_{\textsf{m}}(\vecu, \vecn, \boldsymbol{0}, \boldsymbol{0}) 
- F_{\textsf{m}}(\vecu, \boldsymbol{0}, \boldsymbol{0}, \boldsymbol{0}) \right) \\
&\quad \le \frac{\eta}{e_1} \ell_n \epsilon_n \| \vecu - \vecu^* \|.
\end{align*}
Finally, by Assumption~\ref{as:bound_train_error}, the approximation error satisfies
\[
\left\| \mathcal{F}^{\textsf{NN}}(\vecu, \vecxi, \vectheta) 
- F_{\textsf{ln}}(\vecu, \vecxi, \vectheta) \right\| \le \epsilon^{\textsf{NN}}.
\]
Hence,
\[
2\eta (\vecu - \vecu^*)^\top P 
\left[ \mathcal{F}^{\textsf{NN}}(\vecu, \vecxi, \vectheta) 
- F_{\textsf{ln}}(\vecu, \vecxi, \vectheta) \right] 
\le \frac{\eta}{e_1} \epsilon^{\textsf{NN}} \| \vecu - \vecu^* \|.
\]

Putting all terms together, we get:
\begin{align*}
\dot{V}_1(\vecu) 
&\le -\eta \|\vecu - \vecu^*\|^2 
+ \frac{\eta L}{e_1} \|\vecu - \vecu^*\|^3 \\
&\quad + \frac{\eta}{e_1} \left( \ell_{J_v} E_{J_v} + \ell_{J_i} E_{J_i} 
+ \ell_n \epsilon_n + \epsilon^{\textsf{NN}} \right) \|\vecu - \vecu^*\|.
\end{align*}
We rewrite the inequality by factoring \( \|\vecu - \vecu^*\|^2 \) from the first two terms:
\begin{align*}
\dot{V}_1(\vecu) 
&\le \|\vecu - \vecu^*\|^2 \left( -\eta + \frac{\eta L}{e_1} \|\vecu - \vecu^*\| \right) \\
&\quad + \frac{\eta}{e_1} \left( \ell_{J_v} E_{J_v} 
+ \ell_{J_i} E_{J_i} + \ell_n \epsilon_n + \epsilon^{\textsf{NN}} \right) 
\|\vecu - \vecu^*\|.
\end{align*}
This inequality holds if \( \|\vecu - \vecu^*\| \le \frac{e_1}{L}(1 - s) \), for any \( s \in (s_{\min}, 1] \). Then, the dominant terms yield:
\begin{align*}
\dot{V}_1(\vecu) 
&\le -\eta s \|\vecu - \vecu^*\|^2 \\
&\quad + \frac{\eta}{e_1} 
\left( \ell_{J_v} E_{J_v} + \ell_{J_i} E_{J_i} + \ell_n \epsilon_n + \epsilon^{\textsf{NN}} \right) 
\|\vecu - \vecu^*\|.
\end{align*}
Define \( V_2(\vecu) := \sqrt{V_1(\vecu)} \). Then, using the chain rule,
\[
\dot{V}_2(\vecu) = \frac{\dot{V}_1(\vecu)}{2 V_2(\vecu)}.
\]
Substituting the bound on \( \dot{V}_1(\vecu) \) yields:
\begin{align*}
\dot{V}_2(\vecu) 
& \le -e_1 \eta s V_2(\vecu) \\ 
& + \frac{\eta \sqrt{2e_2}}{2e_1}
\left( \ell_{J_v} E_{J_v} + \ell_{J_i} E_{J_i} 
+ \ell_n \epsilon_n + \epsilon^{\textsf{NN}} \right).
\end{align*}

Let \( b = e_1 \eta s \), and define
\[
a = \frac{\eta \sqrt{2e_2}}{2e_1}
\left( \ell_{J_v} E_{J_v} + \ell_{J_i} E_{J_i} + \ell_n \epsilon_n + \epsilon^{\textsf{NN}} \right).
\]
Then the differential inequality becomes \( \dot{V}_2(\vecu) \le -b V_2(\vecu) + a \), and by Grönwall’s inequality:
\[
V_2(t) \le V_2(t_0) e^{-b(t - t_0)} + \frac{a}{b} \left( 1 - e^{-b(t - t_0)} \right).
\]

Using the bounds \( \|\vecu(t) - \vecu^*\| \le \sqrt{2e_2} V_2(t) \) and \( V_2(t_0) \le \frac{1}{\sqrt{2e_1}} \|\vecu(t_0) - \vecu^*\| \), we obtain:
\begin{align*}
\|\vecu(t) - \vecu^*\| 
&\le \sqrt{\frac{e_2}{e_1}} \|\vecu(t_0) - \vecu^*\| e^{-b(t - t_0)} \\
&\quad + \frac{\sqrt{2e_2}}{b} \cdot a \left( 1 - e^{-b(t - t_0)} \right).
\end{align*}
Define the aggregate error:
\[
\epsilon := \ell_{J_v} E_{J_v} + \ell_{J_i} E_{J_i} + \ell_n \epsilon_n + \epsilon^{\textsf{NN}},
\]
By substituting the definitions of \( a \) and \( b \) and simplifying, we obtain:
\begin{align*}
\frac{\sqrt{2e_2}}{b} \cdot a 
&= \frac{\sqrt{2e_2}}{e_1 \eta s} \cdot \frac{\eta \sqrt{2e_2}}{2e_1} \epsilon 
= \frac{e_2}{e_1^2 s} \epsilon.
\end{align*}
then the final bound becomes:
\begin{align*}
\|\vecu(t) - \vecu^*\| 
&\le \sqrt{\frac{e_2}{e_1}} \|\vecu(t_0) - \vecu^*\| e^{-e_1 \eta s (t - t_0)} \\
&\quad + \frac{e_2}{e_1^2 s} \epsilon \left( 1 - e^{-e_1 \eta s (t - t_0)} \right).
\end{align*}

Evaluating the limit as \( t \to +\infty \) yields the desired local exponential stability result. \hfill \( \triangle \)

\subsection{Proof of Proposition~\ref{prop:invariance}}

The proof leverages Nagumo’s Theorem. Consider the SGF controller \( F_{\textsf{ln}}(\vecu, \vecxi, \vectheta) \) in~\eqref{eq:controller_approximated}; for brevity, let \( \hat{v}_j := |\boldsymbol{\Gamma}_{v,j} \vecu + \bar{v}_j(\vecs_l)| \) denote the linearized voltage magnitude at index \( j \in \mathcal{M} \), and let the measurement noise \( \vecn_v \) satisfy \( \|\vecn_v\|_\infty \le \epsilon_v \), as in Assumption~\ref{as:bound_meas_error}. {Moreover, by Assumption~\ref{as:bound_approx_error}, the error between $|v_j(\vecu;\vecs_l)|$ and its linear approximation  is bounded by \(E_v\). Thus, for every \(j\in\mathcal{M}\),
\begin{equation}\label{eq:Ev_comp_bound}
\hat v_j-E_v \le |v_j(\vecu;\vecs_l)| \le \hat v_j+E_v .
\end{equation}
Next, from Assumption~\ref{as:bound_meas_error}, we have  \(\|\vecn_v\|_\infty \le \epsilon_v\). Applying triangle inequalities to the voltage with measurement noise gives,
\begin{equation}\label{eq:noise_mag_bounds}
|v_j(\vecu;\vecs_l)|-\epsilon_v
\;\le\;
|v_j(\vecu;\vecs_l)|+(\vecn_v)_j
\;\le\;
|v_j(\vecu;\vecs_l)|+\epsilon_v .
\end{equation}
Recall that $\vecV_j=|v_j(\vecu;\vecs_l)|+(\vecn_v)_j$; 
therefore, we have the following:}
\begin{align*}
&-\!\boldsymbol{\Gamma}_{v,j}^\top F_{\textsf{ln}}(\vecu, \vecxi, \vectheta)
\le -\beta(\underline{V}-\vecV_j)\le -\beta \left( \underline{V} - |\hat{v}_j| - \epsilon_v - E_{v} \right), \\
&\boldsymbol{\Gamma}_{v,j}^\top F_{\textsf{ln}}(\vecu, \vecxi, \vectheta)
\le -\beta(\vecV_j-\overline{V}) \le -\beta \left( |\hat{v}_j| - (\overline{V} + \epsilon_v + E_{v}) \right),
\end{align*}
where \( E_{v} \) is defined in Assumption~\ref{as:bound_approx_error}, and it is such that the error between each voltage magnitude and its linear approximation is bounded.   
Now consider the NN-SGF:
\[
F^{\textsf{NN}}(\vecu, \vecxi, \vectheta) = F_{\textsf{ln}}(\vecu, \vecxi, \vectheta) + \Delta F, .
\]
with $\|\Delta F\| \le \epsilon^{\textsf{NN}}$. Next:
\begin{align*}
-\boldsymbol{\Gamma}_{v,j}^\top F^{\textsf{NN}}(\vecu, \vecxi, \vectheta) 
&= -\boldsymbol{\Gamma}_{v,j}^\top F_{\textsf{ln}}(\vecu, \vecxi, \vectheta) 
- \boldsymbol{\Gamma}_{v,j}^\top \Delta F \\
&\le -\beta ( \underline{V} - |\hat{v}_j| - \epsilon_v - E_{v} ) + \|\boldsymbol{\Gamma}_{v,j}\|  \epsilon^{\textsf{NN}} \\
\boldsymbol{\Gamma}_{v,j}^\top F^{\textsf{NN}}(\vecu, \vecxi, \vectheta) 
&= \boldsymbol{\Gamma}_{v,j}^\top F_{\textsf{ln}}(\vecu, \vecxi, \vectheta) 
+ \boldsymbol{\Gamma}_{v,j}^\top \Delta F \\
&\le -\beta ( |\hat{v}_j| - \overline{V} - \epsilon_v - E_{v} ) + \|\boldsymbol{\Gamma}_{v,j}\|  \epsilon^{\textsf{NN}} .
\end{align*}
To ensure the vector field is inward-pointing at the boundary of the voltage constraint set, we then require the linear approximation of the voltage magnitude at node $j \in \mathcal{M}$:
\[
|\hat{v}_j| \in \left[ \underline{V} \!-\! \epsilon_v \!-\! E_{v} \!-\! \frac{\|\boldsymbol{\Gamma}_{v,j}\| \epsilon^{\textsf{NN}}}{\beta}, \ 
\overline{V} \!+\! \epsilon_v \!+\! E_{v} \!+\! \frac{\|\boldsymbol{\Gamma}_{v,j}\| \epsilon^{\textsf{NN}}}{\beta} \right].
\]

A similar argument for the current constraints yields:
\[
|\hat{\imath}_{{j}}| := |\boldsymbol{\Gamma}_{i,{j}} \vecu + \bar{\imath}_{{j}}(\vecs_l)| 
\le \overline{I} + \epsilon_i + E_{i} + \frac{\|\boldsymbol{\Gamma}_{i,{j}}\| \epsilon^{\textsf{NN}}}{\beta}.
\]

Define the inflated constraint bounds:
\begin{align*}
\underline{V}_e &:= \underline{V} - \epsilon_v - E_{v} - \frac{\|\boldsymbol{\Gamma}_{v,j}\| \epsilon^{\textsf{NN}}}{\beta}, \\
\overline{V}_e &:= \overline{V} + \epsilon_v + E_{v} + \frac{\|\boldsymbol{\Gamma}_{v,j}\| \epsilon^{\textsf{NN}}}{\beta}, \\
\overline{I}_e &:= \overline{I} + \epsilon_i + E_{i} + \frac{\|\boldsymbol{\Gamma}_{i,{j}}\| \epsilon^{\textsf{NN}}}{\beta}.
\end{align*}
We define the inflated feasible sets as:
\begin{align*}
\mathcal{S}_{s,\hat{v}} &:= \left\{ \vecu \in \mathcal{C} : \underline{V}_s \leq |\hat{v}_j(\vecu; \vecs_l)| \leq \overline{V}_s, j \in \mathcal{M} \right\}, \\
\mathcal{S}_{s,\hat{\imath}} &:= \left\{ \vecu \in \mathcal{C} : |\hat{i}_j(\vecu; \vecs_l)| \leq \overline{I}_s , j \in \mathcal{L}\right\}, \\
{\mathcal{S}_s} &:= \mathcal{S}_{s,\hat{v}} \cap \mathcal{S}_{s,\hat{\imath}}.
\end{align*}
Since the NN-SGF vector field is strictly inward-pointing, Nagumo’s Theorem implies that \( {\mathcal{S}_s} \) is forward invariant under~\eqref{eq:main_NN}.

We note that the sets \(\mathcal{S}_{s,\hat{v}}\) and \(\mathcal{S}_{s,\hat{\imath}}\) are defined through inequalities involving the linearized voltage magnitude and linearized current magnitude, respectively. These inequalities, however, do not directly translate into constraints expressed in terms of the actual voltage and current magnitudes due to the nonlinear gap between the linearized and actual quantities. 

To obtain a forward-invariant characterization in terms of the actual voltage and current magnitudes, we construct enlarged sets that contain \(\mathcal{S}_s\) and are described by the actual quantities. To this end, observe that
$
|\hat{v}_j| - E_{v,j} \le |v_j| \le |\hat{v}_j| + E_{v,j},
$
and
$
|\hat{\imath}_j| - E_{i,j} \le |\imath_j| \le |\hat{\imath}_j| + E_{i,j}.
$
These bounds imply that 
\(
\mathcal{S}_{s,\hat{v}} \subseteq \mathcal{S}_{e,v}
\)
and
\(
\mathcal{S}_{s,\hat{\imath}} \subseteq \mathcal{S}_{e,i},
\)
which further yields
\(
\mathcal{S}_s \subseteq \mathcal{S}_e,
\)  and the factor \(2\) appearing in the definitions of 
\(\underline V_e\), \(\overline V_e\), and \(\overline I_e\). 

\hfill \(\triangle\)

\end{document}